\documentclass[12pt,a4paper,reqno]{amsart}

\usepackage[english]{babel}
\usepackage{amsmath}
\usepackage{physics}
\usepackage{amssymb}
\usepackage{graphicx}
\usepackage{mathtools}
\usepackage{amsthm}
\usepackage{enumitem}
\usepackage{tikz}
\usepackage{stmaryrd}
\usepackage{mathrsfs}
\usepackage{hyperref}
\usepackage[utf8]{inputenc}
\usepackage{tikz-cd}
\usepackage[T1]{fontenc}
\usepackage{txfonts}
\usepackage{lipsum}
\usepackage[hyperpageref]{backref}
\usetikzlibrary{shapes,arrows,calc}

\numberwithin{equation}{section}
\allowdisplaybreaks
\let\oldpart\part
\def\part{\newpage\oldpart}

\hypersetup{
	linktoc=all,
	colorlinks=true,
	citecolor=magenta,
}

\long\def\eprint#1{\href{https://arxiv.org/abs/#1}{arXiv:#1}}

\makeatletter
\protected\def\Crefnameitem#1#2#3{%
	\@ifnextchar[{\Crefnameitem@first{#1}{#2}{#3}}%
	{\Crefnameitem@first{#1}{#2}{#3}[\textup]}%
}
\def\Crefnameitem@first#1#2#3[#4]{%
	\@ifnextchar[{\Crefnameitem@second{#1}{#2}{#3}{#4}}%
	{\Crefnameitem@second{#1}{#2}{#3}{#4}[(\roman*)]}%
}
\def\Crefnameitem@second#1#2#3#4[#5]{%
	\Crefname{#1enumi}{#2}{#3}%
	\AtBeginEnvironment{#1}{%
		\crefalias{enumi}{#1enumi}%
		\setlist[enumerate,1]{label={#4{#5}.},ref={#5}}%
	}%
}
\makeatother

\theoremstyle{plain}
\newtheorem{thm}{Theorem}[section]
\newtheorem*{theorem*}{Theorem}

\newtheorem{prop}[thm]{Proposition}
\newtheorem{lm}[thm]{Lemma}
\newtheorem{assu}[thm]{Assumption}
\newtheorem{coro}[thm]{Corollary}

\newtheorem{thmalph}{Theorem}

\def\customgenericname{}
\newtheorem{innercustomgeneric}{\customgenericname}
\long\def\newcustomtheorem#1#2{%
	\newenvironment{#1}[1]
	{\def\customgenericname{#2}%
		\def\theinnercustomgeneric{##1}%
		\innercustomgeneric}
	{\endinnercustomgeneric}%
}
\newcustomtheorem{customthm}{Theorem}

\theoremstyle{definition}

\newtheorem{df}[thm]{Definition}
\newtheorem{remark}[thm]{Remark}

\newtheorem{definition}[thm]{Definition}

\long\def\fusion#1#2#3{{\binom{#3}{#1\;#2}}}

\def\N{\mathbb{N}}
\def\Z{\mathbb{Z}}
\def\Q{\mathbb{Q}}

\def\C{\mathbb{C}}

\def\ra{\rightarrow}
\def\ssq{\subseteq}
\def\bs{\backslash}
\def\op{\oplus}
\def\o{\otimes}

\def\<{\langle}
\def\>{\rangle}

\def\ds{\dots}

\def\la{\lambda}

\def\al{\alpha}
\def\si{\sigma}
\def\om{\omega}

\def\b{\beta}

\def\Om{\Omega}

\def\g{\mathfrak{g}}
\def\h{\mathfrak{h}}

\def\sl{\mathfrak{sl}}

\def\vac{\mathbf{1}}
\def\A{\mathsf{A}}

\def\scrU{\mathscr{U}}

\def\L{\mathsf{L}}

\def\Mod{\mathsf{Mod}}

\def\Hom{\mathrm{Hom}}

\def\End{\mathrm{End}}

\def\Id{\mathrm{Id}}
\def\wt{\mathrm{wt}}

\def\ker{\operatorname{Ker}}

\def\CC{\mathcal{C}}

\def\bt{\boxtimes}

\newcommand{\Oord}{\mathcal{O}_{\ell,\mathrm{Ord}}}
\newcommand{\Lh}[1]{\hat{L}_\ell(#1)}
\theoremstyle{plain}
\newcustomtheorem{customcoro}{Corollary}
\theoremstyle{definition}

\begin{document}
	
	\title{A Kazhdan--Lusztig correspondence at admissible levels in type $A$}
	
	\author{Jianqi Liu}
	\address{Department of Mathematics, New Uzbekistan University, Tashkent, Uzbekistan}
	\email{j.liu@newuu.uz}
	
	\author{Xingting Wang}
	\address{Department of Mathematics, Louisiana State University, Baton Rouge, LA}
	\email{xingtingwang@lsu.edu}

	\begin{abstract}
Let $\ell=-N+p/v$ be an admissible level for $\hat{\sl}_N$, and let $\Oord$ be the category of ordinary modules over the simple affine vertex operator algebra  $L_\ell(\sl_N)$, with the Huang--Lepowsky--Zhang tensor product. In this paper, we prove that $\Oord$ is ribbon equivalent to the semisimplification of the category of tilting modules for Lusztig's quantum group $U_q(\sl_N)$ at the root of unity $q=e^{\pi iv/p}$. The equivalence sends each irreducible module $\Lh{\la}$ to the simple object $L_q(\la)$ with the same highest weight. This proves the correspondence conjectured by Creutzig--Huang--Yang in type~$A$ with the choice $q^{1/N}=e^{\pi iv/(Np)}$ for the quantum group braiding. For $p>N$, a recognition theorem for ribbon fusion categories with fusion rings isomorphic to that of $L_{p-N}(\sl_N)$ reduces the proof to comparing the twist and self-braiding of $\Lh{\om_1}$ with those of $L_q(\om_1)$, which is achieved by computing the twist and the braiding eigenvalues from the bottom-degree Zhu algebra modules via the $\odot$-product. 
	\end{abstract}
	
	\maketitle

	\section{Introduction}

The study of relations between braided tensor categories arising from different quantum algebras has long been a central topic in representation theory and mathematical physics. 
Two of the best-known examples of quantum algebras arising in quantum conformal field theory are vertex operator algebras (VOAs) \cite{Bor86,FLM88} and quantum groups \cite{Dr85}. Both give rise to braided tensor categories of modules \cite{HL94,H05,RT90} and, under suitable conditions, to modular tensor categories \cite{H08,RT91,TW93}.

For a complex simple Lie algebra $\g$, Drinfeld proved a braided equivalence over $\C[\![\hbar]\!]$ between the category defined by the Knizhnik--Zamolodchikov (KZ) local systems on configuration spaces and the corresponding category of $U_\hbar(\g)$-modules \cite{Dr90}.
For a simply laced affine Lie algebra $\hat{\g}$ at level $k$ with $k+h^\vee\notin\Q_{\geq0}$, Kazhdan and Lusztig constructed a braided tensor category of finite-length $\hat{\g}$-modules using coinvariants with the KZ connection \cite{KL93a,KL93b,TUY89}. They also constructed a braided tensor equivalence with the corresponding category of quantum group modules at generic levels  \cite{KL94a,KL94b}.
Finkelberg later established a positive integral level version of this equivalence between the fusion category of integrable highest-weight modules over the affine Lie algebra \cite{Kac90} and the semisimplification of the category of quantum group modules at a root of unity \cite{F96,F13}.

The representation theory of affine Lie algebras $\hat{\g}$ is closely related to that of the corresponding affine VOAs $L_k(\g)$ \cite{FZ92,DLM97}. Using intertwining operators \cite{FHL} and the analytic properties of their products and iterates, Huang and Lepowsky developed a tensor product theory for modules over strongly rational VOAs \cite{HL94,HL95a,HL95b,H05}, which include the affine VOA $L_k(\g)$ at positive integral level $k$. Huang, Lepowsky and Zhang later extended this theory to suitable categories of
generalized modules using logarithmic intertwining operators \cite{HLZ}.
Using this tensor product theory, Creutzig, Huang and Yang \cite{CHY18} proved that the category of ordinary modules $\Oord$ over the simple affine VOA $L_\ell(\g)$ at admissible level $\ell$ in the sense of Kac and Wakimoto \cite{KaW88}, which is in general neither rational nor $C_2$-cofinite \cite{A16,AvE23}, is a braided tensor category, and conjectured that $\Oord$ is braided equivalent to the semisimplification of the category of tilting modules for the corresponding quantum group $U_q(\g)$ at a root of unity \cite[Conjecture~1.3]{CHY18}, generalizing the braided equivalence of Kazhdan--Lusztig and Finkelberg to admissible levels.

In this paper, we establish the correspondence proposed by Creutzig, Huang and Yang for $\sl_N$ at every admissible level $\ell$, with $q=e^{\pi i/(\ell+N)}$ and $q^{1/N}=e^{\pi i/(N(\ell+N))}$ for the quantum group braiding.
Our proof uses the recent identification of the bottom degree of the fusion tensor product $M^1\bt M^2=M^1\bt_{P(1)}M^2$ with the contracted tensor product $M^1\odot M^2$ via restricted conformal blocks \cite{GLZ25,GLZ26,Liu26,CL26}, together with the reconstruction theorem of Kazhdan and Wenzl \cite{KW93} in the nonpointed case.

	\subsection{The main theorem and consequences}
We state our main theorem first. Let $N\geq2$. 
Following \cite{A16,CHY18}, we write the admissible level of $\hat{\sl}_N$ as
	\[
		\ell=-N+\frac{p}{v},\qquad p,v\in\Z_{\geq1},\quad \gcd(p,v)=1,\quad p\geq N.
	\]
	Let $k=p-N$. The irreducible ordinary $L_\ell(\sl_N)$-modules are the modules $\Lh{\la}$ whose bottom degrees are the irreducible $\sl_N$-modules $L(\la)$, with
	\[
		\la\in P^k_+
		=\left\{\sum_{i=1}^{N-1}a_i\om_i:a_i\in\Z_{\geq0},\ \sum_{i=1}^{N-1}a_i\leq k\right\},
	\]
where $\om_1,\ldots,\om_{N-1}$ are the fundamental weights of $\sl_N$. In particular, at all admissible levels with the same numerator $p$, the irreducible ordinary modules have the same labels and the same fusion rules, namely those of $L_k(\sl_N)$ \cite{A16,Cre19}.

Let $q=e^{\pi iv/p}$. Denote by $\CC(\sl_N,q)$ the quotient of the category of type~$1$ tilting modules for $U_q(\sl_N)$ by the ideal of negligible morphisms \cite{AP95,Saw06}. Its simple objects are the images $L_q(\la)$ of the Weyl modules with $\la\in P^k_+$. 
Choose an $N$-th root $y$ of $q$. 
For tilting modules $U,V$, the Cartan factor of the $R$-matrix acts on $U_\al\o_\C V_\beta$ by $y^{N(\al\mid\beta)}$. The $R$-matrix followed by the flip induces the braiding $c^y$ on $\CC(\sl_N,q)$ as in \eqref{eq:QG-braiding}. The following is our main theorem. 

	\begin{thmalph}[Theorem~\ref{thm:main}]\label{thm:intro-ribbon-equivalence}
		Let $q=e^{\pi iv/p}$ and $y=e^{\pi iv/(Np)}$. Then there is a ribbon equivalence
		\[
			F:\Oord\xrightarrow{\simeq}(\CC(\sl_N,q),c^y)
		\]
		such that $F(\Lh{\la})\cong L_q(\la)$ for every $\la\in P^k_+$.
	\end{thmalph}

Note that for $k\geq1$, different choices of the $N$-th root $y$ of $q$ give different self-braidings of $L_q(\om_1)$ by Lemma~\ref{lm:QG-braiding}. 
The choice of $y$ in Theorem~\ref{thm:intro-ribbon-equivalence} is determined by comparing the self-braiding eigenvalues of $L_q(\om_1)$ and $\Lh{\om_1}$ on the corresponding summands.

When $v=1$, we have $\ell=k$ and $q=e^{\pi i/(k+N)}$. Therefore, for $k\geq1$, Theorem~\ref{thm:intro-ribbon-equivalence} recovers the type~$A$ positive integral level case of the correspondence in \cite{KL93a,KL93b,KL94a,KL94b,F96,F13}. For general $v$, we can also use the equivalence to determine when $\Oord$ is modular. More precisely, Corollary~\ref{coro:modularity} gives
	\[
		\Oord\text{ is modular}\quad\Longleftrightarrow\quad\gcd(v,N)=1.
	\]
The implication from right to left was proved by Creutzig \cite[Theorem~3]{Cre19}. 
For $p>N$, the converse follows by applying the type~$A$ modularity criterion for the quantum group category to the braiding $c^y$ in Theorem~\ref{thm:intro-ribbon-equivalence} \cite{Bru00,Row06}. When $p=N$, both categories are equivalent to $\mathrm{Vec}$ and $\gcd(v,N)=1$ by admissibility.

	\subsection{Main ideas of the proof}

The categories $\Oord$ and $(\CC(\sl_N,q),c^y)$ have the same fusion rules \cite{AP95,Cre19}, and for $k\geq1$ their simple objects occur in tensor powers of $X=\Lh{\om_1}\in \Oord$ and $Y=L_q(\om_1)\in (\CC(\sl_N,q),c^y)$, respectively. We construct the ribbon equivalence by comparing the twists of $X$ and $Y$ and the eigenvalues of $c_{X,X}$ and $c^y_{Y,Y}$ on the corresponding simple summands of $X\bt X$ and $Y\bt Y$.

In order to determine the braiding eigenvalues for $\Oord$, we use the explicit braiding formula on the bottom degree of the Huang--Lepowsky tensor product $M^1\bt_{P(1)} M^2$ \cite{HL95a}, which is identified with $M^1\odot M^2$ \cite{Liu26,AL26}, a left module over the Zhu algebra $\A\simeq U(\sl_N)/J_\ell$ \cite{FZ92,Z96}. Note that the category $\Oord$ is semisimple \cite{A16}, and it is equivalent to the category of finite-dimensional modules over $\A$ via the generalized Verma module functor $\Phi^\L$ \cite{DLM98,DGK25,AL26}.
 For irreducible modules $\Lh{\la},\Lh{\mu}\in \Oord$, let $\varphi_{\la,\mu}:L(\la)\o_\C L(\mu)\twoheadrightarrow\Lh{\la}\odot\Lh{\mu}$ be the map $u\o w\mapsto u\odot w$. Under this identification, we prove in Theorem~\ref{prop:adm-braiding} that the braiding on bottom degrees is given by
\begin{equation}\label{eq:intro-adm-braiding}
	\mathcal R_{\Lh{\la},\Lh{\mu}}(u\odot w)
	=\varphi_{\mu,\la}\left(\tau\exp\left(\frac{\pi i}{\ell+h^\vee}\Om_{\la,\mu}\right)(u\o w)\right),
\end{equation}
where $\tau(u\o w)=w\o u$ and $\Om_{\la,\mu}(u\o w)=\sum_a x_a u\o x^a w$ for dual bases $\{x_a\}$ and $\{x^a\}$ of $\g$. For the formula at positive integral levels, see \cite[Theorem~5.4]{CL26}. In particular, on a multiplicity-one summand $\Lh{\nu}$ of $\Lh{\la}\bt\Lh{\la}$, the exponential acts by $e^{\pi i(h_\nu-2h_\la)}$, and the flip $\tau$ determines the sign of the braiding eigenvalue.

In Section~\ref{sec:recognition}, we prove a general theorem that compares a ribbon fusion category $\mathcal B$ with the fusion rules of $L_k(\sl_N)$ to the quantum group side $(\CC(\sl_N,q),c^y)$. For $k\geq2$, the reconstruction theorem of Kazhdan and Wenzl \cite{KW93} identifies the monoidal category up to a twist of the
associator by a $3$-cocycle on $\Z/N$. Let $X\in \mathcal B$ and $Y=L_q(\om_1)\in (\CC(\sl_N,q),c^y)$ be the respective simple objects corresponding to $\om_1$. 
For $k\geq2$, suppose that their self-braidings on $X$ and $Y$ have the same eigenvalues $\lambda_+$ and $\lambda_-$ on the summands corresponding to $2\om_1$ and $\om_2$, respectively. Then both categories have Hecke parameter $Q=-\lambda_-/\lambda_+$ and determinant scalar $(-\lambda_+)^{-N}$ for this choice of $Q$. 
By \cite[Theorem~4.1]{KW93}, the common parameter $Q$ identifies the endomorphism algebras of $X^{\bt m}$ and $Y^{\bt m}$ as the same quotient 
of $H_m(Q)$ for each $m\geq2$. 
In Proposition~\ref{thm:KW-comparison}, we use \cite[Proposition~4.6]{TW05} and equality of the determinant scalars to extend $\End_{\mathcal B}(X^{\bt m})\cong\End_{\CC(\sl_N,q)}(Y^{\bt m})$ to a label-preserving tensor
equivalence $F$ from $\mathcal B$ to $\CC(\sl_N,q)$. In this way, 
we avoid a direct comparison of the associators, which on the VOA side would involve products and iterates of intertwining operators \cite{H05,HLZ,MS89}.
Under the identification $F(X\bt X)\cong Y\bt Y\cong L_q(2\om_1)\op L_q(\om_2)$, both $F(c_{X,X})$ and $c^y_{Y,Y}$ act by $\lambda_+$ on $L_q(2\om_1)$ and by $\lambda_-$ on $L_q(\om_2)$. Hence $F(c_{X,X})=c^y_{Y,Y}$. 
Applying the hexagon identities, we see that $F$ preserves the braiding on tensor powers and direct summands  of $X$. 

The following criterion in Theorem~\ref{thm:recognition} also includes the pointed case $k=1$: 
 Suppose that $\mathcal B$ is a ribbon fusion category whose fusion ring is identified with $\mathfrak F(L_k(\sl_N))$ as a based ring. If the twist of $X=X_{\om_1}$ is $y^{N^2-1}\Id_X$ and $c_{X,X}$ acts by $-y^{-N-1}$ on $X_{\om_2}$ and, for $k\geq2$, by $y^{N-1}$ on $X_{2\om_1}$, then there is a ribbon equivalence
\[
	F:\mathcal B\xrightarrow{\simeq}(\CC(\sl_N,q),c^y),\qquad
	F(X_\la)\cong L_q(\la)\quad\text{for every}\ \la\in P^k_+.
\]
Finally, in Section~\ref{sec:main-section}, we apply the criteria  Theorem~\ref{thm:recognition} to $\mathcal B=\Oord$. By Lemma~\ref{lm:adm} and \eqref{eq:QG-fusion-rule}, both categories have the fusion rules of $L_k(\sl_N)$ \cite{Bea96,V88}. For $k\geq1$, let $X=\Lh{\om_1}$. In Proposition~\ref{prop:adm-invariants}, we use our explicit braiding formula \eqref{eq:intro-adm-braiding} on the bottom degree of $X\bt X$ to compute the eigenvalues of $c_{X,X}$. 
These eigenvalues and the twist of $X$ in \eqref{eq:adm-vector-twist} show that $\Oord$ satisfies the hypotheses of Theorem~\ref{thm:recognition} for $k\geq1$.


	
	\section{Preliminaries}\label{sec:prelim}
Throughout this paper, all categories are $\C$-linear, and all vector spaces are over $\C$.
	
	\subsection{Affine VOAs at admissible levels}\label{subsec:adm}
Let $\g$ be a finite-dimensional simple Lie algebra.
Fix a Cartan subalgebra $\h$ and a set of positive roots $\Delta_+$. Let $P_+$ be the set of dominant integral weights, let $\theta$ be the highest root, and let $\rho=\frac12\sum_{\al\in\Delta_+}\al$. 
We normalize the invariant bilinear form $(\cdot\mid\cdot)$ so that $(\theta\mid \theta)=2$. Let $w_0$ be the longest element of the Weyl group, so that $w_0(\Delta_+)=-\Delta_+$. For $\la\in P_+$, let $L(\la)$ be the finite-dimensional irreducible $\g$-module of highest weight $\la$; the Casimir element acts on it by $(\la\mid\la+2\rho)$. Denote by $h$, $h^\vee$ and $r^\vee$ the Coxeter number, dual Coxeter number and lacing number of $\g$, respectively.

	The affine Lie algebra is defined by $\hat{\g}=\g\o\C[t,t^{-1}]\op\C K$, with Lie bracket \cite[Chapter~7]{Kac90}
	\[
	[x\o t^m,y\o t^n]=[x,y]\o t^{m+n}+m\delta_{m+n,0}(x\mid y)K,\qquad [K,\hat{\g}]=0.
	\]
	We write $x(n)=x\o t^n$ for $x\in\g$ and $n\in\Z$. For $\ell\in\C$, the induced module
	\[
	V^\ell(\g)=U(\hat{\g})\o_{U(\g[t]\op\C K)}\C\vac
	\]
	is the {\bf universal affine vertex algebra}, where $\g[t]$ acts trivially on $\C\vac$ and $K$ acts as $\ell\cdot\Id$. Its generating fields are $Y(x(-1)\vac,z)=\sum_{n\in\Z}x(n)z^{-n-1}$. For $\ell\neq-h^\vee$, we equip $V^\ell(\g)$ with its Sugawara conformal vector and denote its unique simple quotient by $L_\ell(\g)$ \cite{FZ92}.
	
	Following \cite{KaW88,KaW89}, a level $\ell$ is {\bf admissible} if and only if
	\begin{equation}\label{eq:kappa}
		\kappa=\ell+h^\vee=\frac{p}{v},\qquad p,v\in\Z_{\geq1},\quad \gcd(p,v)=1,
		\qquad p\geq
		\begin{cases}
			h^\vee,&\gcd(r^\vee,v)=1,\\
			h,&\gcd(r^\vee,v)=r^\vee;
		\end{cases}
	\end{equation}
	From now on, assume that $\ell$ is admissible and $V=L_\ell(\g)$.
	
	Choose dual bases $\{x_a\}$ and $\{x^a\}$ of $\g$ with respect to $(\cdot\mid\cdot)$. By the Segal--Sugawara construction, the Virasoro element of $L_\ell(\g)$ is
	\begin{equation}\label{eq:sugawara}
		\om=\frac{1}{2\kappa}\sum_a x^a(-1)x_a(-1)\vac.
	\end{equation} We write $Y(\om,z)=\sum_{n\in \Z} \om_{(n)}z^{-n-1}=\sum_{n\in\Z}L(n)z^{-n-2}$.
	Given a dominant weight $\la\in P_+$, denote by $\Lh{\la}$ the irreducible quotient $\hat{\g}$-module of $U(\hat{\g})\o_{U(\g[t]\op\C K)}L(\la)$, where $\g\o t\C[t]$ acts trivially on $L(\la)$ and $K$ acts as $\ell\cdot\Id$. Then $	\Lh{\la}$ is an $\N$-gradable module over the VOA $V^\ell(\g)$: 
	\[
	\Lh{\la}=\bigoplus_{n\geq0}\Lh{\la}(n),\qquad \Lh{\la}(0)=L(\la),
	\]
	where $\Lh{\la}(n)$ is the $\om_{(1)}=L(0)$ eigenspace of eigenvalue $h_\la+n$. The conformal weight $h_\la$ is
	\begin{equation}\label{eq:conformal-weight}
		h_\la=\frac{(\la\mid\la+2\rho)}{2\kappa}=\frac{v(\la\mid\la+2\rho)}{2p}.
	\end{equation}
	In particular, as an irreducible $\hat{\g}$-module, $L_\ell(\g)=\Lh{0}$.
	
	The set of ordinary admissible weights at level $\ell$ is
	\begin{equation}\label{eq:PrZ}
		\mathrm{Pr}^{\ell}_{\Z}=\{\la\in P_+:\Lh{\la}\text{ is an }L_\ell(\g)\text{-module}\};
	\end{equation}
	see \cite[Section~2 and Theorem~3.1]{AvE23}. Let $\Oord$ denote the category of ordinary $L_\ell(\g)$-modules in the  category $\mathcal O$ of affine Lie algebra $\hat{\g}$.
	
	\begin{lm}[\cite{A16}]\label{lm:adm-semisimple}
		The category $\Oord$ is a finite semisimple category, with irreducible objects $\{\Lh{\la}:\la\in\mathrm{Pr}^{\ell}_{\Z}\}$. Every object of $\Oord$ is a finite direct sum of these modules.
	\end{lm}
	
	Let $A(V)=V/O(V)$ be the Zhu algebra of a VOA $V$ \cite{Z96}, and let $[a]$ denote the image of $a\in V$ in $A(V)$. Recall that $x\mapsto[x(-1)\vac]$ induces an isomorphism $U(\g)\cong A(V^\ell(\g))$, for $x\in \g$ \cite{FZ92}. Passing to the simple quotient VOA $L_\ell(\g)$, we have
	\[
	\A:=A(L_\ell(\g))\cong U(\g)/J_\ell,\qquad [x(-1)\vac]\mapsto x+J_\ell,\quad x\in\g,
	\]
	where $J_\ell$ is the kernel of the natural surjective homomorphism
	\[
	U(\g)\cong A(V^\ell(\g))\ra A(L_\ell(\g)).
	\]
	In particular, every left $\A$-module is a $\g$-module.
	
	For $M\in\Oord$, recall that the bottom degree subspace is a  left $\A$-module
	\[
	\Om(M)=\{v\in M:a_{(n)}v=0,\ \forall\ a\in V\text{ and }n\in\Z\text{ with }\wt a-n-1<0\},
	\]
	with $[a]$ acting as $o(a)=a_{(\wt a-1)}$ for homogeneous $a\in L_\ell(\g)$. In particular, $[x(-1)\vac]$ acts as $x(0)$ for $x\in\g$. Note that for any $\la\in \Pr_\Z^\ell$, we have $\Om(\Lh{\la})= \Lh{\la}(0)=L(\la)$, and $[\om]\in \A$ acts on this space by $h_\la\cdot \Id$. In particular, we have a well-defined functor 
	\begin{equation}
		\Om: \Oord\to \mathsf{Mod}_{\mathrm{fd}}(\A),
	\end{equation}
	where  $\mathsf{Mod}_{\mathrm{fd}}(\A)$ is the category of finite-dimensional left $\A$-modules. Every such module belongs to the BGG category $\mathcal O$ of $\g$, since it is completely reducible as a finite-dimensional $\g$-module. This is the category denoted by $\mathcal O_{\ell,\mathrm{Ord}}^{\mathsf{fin}}$ in \cite[Definition~3.1(4)]{AL26}.
	
	On the other hand, given $S\in \Mod_{\mathrm{fd}}(\A)$, the generalized Verma module \cite{DLM98,DGK25} 
	\[
	\Phi^\L(S)=(\scrU/\mathsf N^1_{\mathsf{L}}\scrU)\otimes_{\scrU_0} S
	\]
	is an $\N$-gradable $L_\ell(\g)$-module, where $\scrU=\scrU(L_\ell(\g))$ is the universal enveloping algebra of the VOA $L_\ell(\g)$, and $\mathsf N^1_{\mathsf{L}}\scrU$ is the first neighborhood basis \cite{FZ92,FBZ04}. 
	\begin{lm}[{\cite[Theorem~3.2 and Corollary~3.3]{AL26}}]\label{lm:adm-Zhu-equivalence}
		For any $S\in \Mod_{\mathrm{fd}}(\A)$, the generalized Verma module $\Phi^\L(S)$ is in $\Oord.$ Moreover, the functor $\Phi^\L$ and the bottom-degree functor $\Om$ form an adjoint equivalence
		\begin{equation}\label{eq:adm-Zhu-equivalence}
			(	\Phi^\L\dashv \Om):\mathsf{Mod}_{\mathrm{fd}}(\A)\rightleftarrows\Oord.
		\end{equation}
		Both categories are semisimple, and their irreducible objects correspond by
		\begin{equation}\label{eq:adm-Zhu-simples}
			\Phi^\L(L(\la))\cong\Lh{\la},\qquad
			\Om(\Lh{\la})=L(\la),\qquad \la\in \mathrm{Pr}^{\ell}_{\Z}.
		\end{equation}
	\end{lm}
	In particular, every finite-dimensional left $\A$-module is a finite direct sum of the modules $L(\la)$ with $\la\in \mathrm{Pr}^{\ell}_{\Z}$. Moreover, the natural maps $S\ra\Om(\Phi^\L(S))$ and $\Phi^\L(\Om(M))\ra M$ are isomorphisms for $S\in\mathsf{Mod}_{\mathrm{fd}}(\A)$ and $M\in\Oord$.

	\subsection{Huang--Lepowsky tensor products in \texorpdfstring{$\Oord$}{O(ell,Ord)}}\label{subsec:HL}
	We recall Huang--Lepowsky's $P(w)$-tensor product for modules in the category $\Oord$; see \cite{HL94,HL95a,HL95b,HLZ}. 
	
	By Lemma~\ref{lm:adm-semisimple}, every $M\in\Oord$ is a finite direct sum of irreducible modules. Hence $M$ is $\N$-gradable, with $M=\bigoplus_{n\geq0}M(n)$ and $M(0)=\Om(M)$. Denote the completion of $M$ and its bottom-degree projection by
	\begin{equation}\label{eq:HL-bottom-projection}
		\overline M=\prod_{n\geq0}M(n),\qquad
		\pi_0:\overline M\ra M(0)=\Om(M).
	\end{equation}
	For a $V$-module homomorphism $f:M\ra N$, denote its extension to the completions by $\bar f:\overline{M}\to \overline{N}$. Since $f$ preserves the relative grading, $L(-1) M(n)\ssq M(n+1)$, and $L(0)$ acts on $\Om(M)$ as $[\om]$, then for any $y\in \C$ we have
	\begin{equation}\label{eq:HL-projection-identities}
		\pi_0\bar f=\Om(f)\pi_0,\qquad
		\pi_0e^{yL(0)}=e^{y[\om]}\pi_0,\qquad
		\pi_0e^{yL(-1)}=\pi_0.
	\end{equation}
	Note that the expression $e^{yL(-1)}v$ is well-defined for $v\in \overline M$, since only finitely many summands contribute to each homogeneous component $M(n)$. See also \cite[Section~3.3.3]{CL26}.
	
	For $M^1,M^2,M^3\in\Oord$, let $I\fusion{M^1}{M^2}{M^3}$ be the space of formal intertwining operators of the indicated type \cite{FHL}. For $w\in\C^\times$, fixing a value of $\log w$ identifies intertwining operators with $P(w)$-intertwining maps by evaluation at $e^{\log w}$ \cite{HL95a,HL95b,HLZ}.
	
	\begin{df}[{\cite{HL95a,HL95b,HLZ}}]\label{df:HL-tensor-product}
		Let $M^1,M^2\in\Oord$ and $w\in\C^\times$. A {\bf $P(w)$-tensor product} of $M^1$ and $M^2$ consists of a module $M^1\bt_{P(w)}M^2\in\Oord$ and a canonical $P(w)$-intertwining map
		\begin{equation}\label{eq:HL-Pz-map}
			\bt_{P(w)}:M^1\o_\C M^2\ra\overline{M^1\bt_{P(w)}M^2},
			\qquad v^1\o v^2\mapsto v^1\bt_{P(w)}v^2,
		\end{equation}
		satisfying the following universal property: For a fixed value of $\log w$, let
		$\mathcal Y_{12}\in
		I\fusion{M^1}{M^2}{M^1\bt_{P(w)}M^2}$
		be the unique intertwining operator satisfying
		\[
		\mathcal Y_{12}(v^1,e^{\log w})v^2
		=v^1\bt_{P(w)}v^2,
		\qquad v^i\in M^i.
		\]
		Then there is a unique natural isomorphism given by composition:
		\begin{equation}\label{eq:adm-intertwiner-universal}
			\Hom_{\Oord}(M^1\bt_{P(w)}M^2,M^3)\xrightarrow{\simeq}
			I\fusion{M^1}{M^2}{M^3},\qquad f\mapsto f\circ\mathcal Y_{12}.
		\end{equation}
	\end{df}
	
	We take $w=1$ and $\log(1)=0$, and write $\bt=\bt_{P(1)}$. Then the corresponding intertwining operator $\mathcal Y_{12}$ satisfies
	\[
	\mathcal Y_{12}(v^1,e^0)v^2=v^1\bt v^2\in\overline{M^1\bt M^2}.
	\]
	
	Recall that the braiding $c_{M^1,M^2}$ in the category $\Oord$ is defined via the commutativity isomorphism of the space of intertwining operators followed by parallel transport from $-1$ to $1$ through the upper half-plane. With $\log(1)=0$, the corresponding value of the logarithm at $-1$ is $\pi i$. Let $\mathcal Y_{21}$ be the canonical intertwining operator with target $M^2\bt M^1$. Then
	\begin{equation}\label{eq:HL-braiding}
		\overline{c_{M^1,M^2}}(v^1\bt v^2)
		=e^{L(-1)}\mathcal Y_{21}(v^2,e^{\pi i})v^1,
	\end{equation}
	where the evaluation at $e^{\pi i}$ uses this value of the logarithm; see \cite[Part~VIII, Section~12.2]{HLZ}.
	
	We will use the following result.
	\begin{lm}[{\cite[Theorems~4.3 and~6.6]{CHY18}}]\label{lm:HL-tensor}
		The category $\Oord$ is a braided tensor category with the Huang--Lepowsky tensor product $\bt=\bt_{P(1)}$, tensor unit $L_\ell(\g)$, and braiding \eqref{eq:HL-braiding}.
	\end{lm}

	\subsection{Ribbon fusion categories}\label{subsec:ribbon}
	We refer to \cite{EGNO,BK01,Tur94} for the definitions of fusion categories, braided fusion categories and ribbon categories. 
	
Let $\mathcal B$ be a fusion category with tensor product $\bt$ and unit object $\vac$. Let $X_0=\vac,X_1,\ldots,X_r$ be representatives of the isomorphism classes of simple objects of $\mathcal B$. The {\bf fusion ring} $\mathfrak F(\mathcal B)$ is the Grothendieck ring of $\mathcal B$. As an abelian group, it is free with distinguished basis $\{[X_i]:0\leq i\leq r\}$, and its multiplication is given by
	\begin{equation}\label{eq:fusion-ring-product}
		[X_i]\cdot[X_j]=\sum_{s=0}^r N_{i,j}^s[X_s],\qquad
		N_{i,j}^s=\dim\Hom_{\mathcal B}(X_i\bt X_j,X_s).
	\end{equation}
	Its unit is $[\vac]$, and duality induces the involution $[X_i]\mapsto[X_i^*]$. If $\mathcal B$ is braided, then $\mathfrak F(\mathcal B)$ is commutative. We say that $[X_s]$ {\em occurs} in a product of basis elements if its coefficient in the distinguished basis is positive. An isomorphism of fusion rings is an isomorphism {\bf of based rings} if it maps one distinguished basis bijectively onto the other.
	
	Let $\mathcal B$ be a braided tensor category with braiding $c$. Recall that a twist $\theta$ is a natural automorphism of the identity functor such that $\theta_\vac=\Id_\vac$ and
	\begin{equation}\label{eq:balancing}
		\theta_{X\bt Y}=c_{Y,X}\circ c_{X,Y}\circ(\theta_X\bt\theta_Y),\qquad X,Y\in\mathcal B.
	\end{equation}
	If $\mathcal B$ is rigid, the twist is a {\bf ribbon twist} if $\theta_{X^*}=(\theta_X)^*$ for every $X\in\mathcal B$. For a simple object $X$ with $\End_{\mathcal B}(X)=\C$, we also write $\theta_X$ for the scalar by which the twist acts on $X$.
	
	\begin{df}\label{df:braiding-eigenvalue}
		Let $\mathcal B$ be a semisimple braided tensor category whose simple objects have endomorphism algebra $\C$. Let $X,Z$ be simple objects of $\mathcal B$, and assume that $Z$ occurs with multiplicity one in $X\bt X$. For any nonzero morphism $\iota:Z\ra X\bt X$, there is a unique scalar $\lambda_Z\in\C^\times$ such that
		\begin{equation}\label{eq:braiding-eigenvalue}
			c_{X,X}\circ\iota=\lambda_Z\iota.
		\end{equation}
		Since $\dim\Hom_{\mathcal B}(Z,X\bt X)=1$, the scalar is independent of the choice of $\iota$. We call $\lambda_Z$ the {\bf eigenvalue of $c_{X,X}$ on $Z$}.
	\end{df}
	
	\begin{lm}\label{lm:ribbon-structures}
		Let $(\mathcal B,c)$ be a braided fusion category and let $X\in\mathcal B$. Suppose that every simple object of $\mathcal B$ occurs in $X^{\bt n}$ for some $n\geq0$. If $\theta$ and $\theta'$ are twists on $(\mathcal B,c)$ and $\theta_X=\theta'_X$, then $\theta=\theta'$.
	\end{lm}
	\begin{proof}
		By \eqref{eq:balancing}, if $\theta_{X^{\bt n}}=\theta'_{X^{\bt n}}$, then
		\[
		\theta_{X^{\bt(n+1)}}=c_{X,X^{\bt n}}\circ c_{X^{\bt n},X}\circ
		(\theta_{X^{\bt n}}\bt\theta_X)=\theta'_{X^{\bt(n+1)}}.
		\]
		Since $\theta_\vac=\theta'_\vac=\Id_\vac$, the twists agree on $X^{\bt n}$ for every $n\geq0$. Let $Z$ be a simple object and choose an embedding $\iota:Z\hookrightarrow X^{\bt n}$. By naturality, we have
		\[
		\iota\circ\theta_Z=\theta_{X^{\bt n}}\circ\iota
		=\theta'_{X^{\bt n}}\circ\iota=\iota\circ\theta'_Z.
		\]
		Hence $\theta_Z=\theta'_Z$. Since $\mathcal B$ is semisimple, the two twists agree on every object.
	\end{proof}
	
	For a ribbon fusion category $\mathcal B$, the {\bf M\"uger center} is the full subcategory of objects $Z$ such that $c_{Y,Z}\circ c_{Z,Y}=\Id_{Z\bt Y}$ for every $Y\in\mathcal B$. The category $\mathcal B$ is {\bf modular} if its M\"uger center consists of finite direct sums of the unit object \cite{Mu03,DGNO10}.
	
	\begin{lm}[{\cite[Corollary~7.3]{Cre19}}]\label{lm:adm-ribbon-ADE}
		Let $\g$ be a simple Lie algebra of type $A$, $D$ or $E$, and let $\ell$ be an admissible level. Then $\Oord$ is a ribbon fusion category with twist $\theta_M=e^{2\pi iL(0)}$ for $M\in\Oord$. In particular, $\Oord$ is rigid.
	\end{lm}
	
	
	\section{Braiding of ordinary modules at admissible levels}\label{sec:admissible}
	In this section, we give an explicit formula for the braiding on the bottom degrees of ordinary modules for an admissible affine VOA and determine the braiding eigenvalues on simple summands of multiplicity one.
	
	Let $\g$ be a finite-dimensional complex simple Lie algebra and let $\ell$ be an admissible level. We use the notation of Subsections~\ref{subsec:adm} and~\ref{subsec:HL}; in particular, $V=L_\ell(\g)$ and $\kappa=\ell+h^\vee$.
	
	\subsection{Bottom degrees of tensor products in $\Oord$}

	We first recall the contracted tensor product of VOA modules introduced by the first-named author in \cite{Liu26}.
	
	\begin{df}[{\cite[Definition~4.3]{Liu26}}]\label{df:contracted-tensor-product}
		Let $M^1,M^2\in\Oord$. Define
		\[M^1\odot M^2:=\frac{M^1\o_\C M^2}
		{\mathcal{L}_{\mathbb P^1\bs\{1,0\}}(V)_{<0}.(M^1\o_\C M^2)},\]
		where $\mathcal{L}_{\mathbb P^1\bs\{1,0\}}(V)_{<0}$ is the subalgebra of the three-point chiral Lie algebra $\mathcal{L}_{\mathbb P^1\bs\{\infty,1,0\}}(V)$ \cite{NT05,DGT24} defined by a condition on the order of poles at $\infty$. It is spanned by the classes
		\[X=a\o\frac{z^n}{(z-1)^m} (dz)^{1-\wt a},\qquad a\in V,\quad m,n\in\Z,\quad n-m<\wt a-1,	\]
		where $a\in V$ is homogeneous and $z$ is the coordinate on $\mathbb P^1\bs \{ \infty,1,0\}$. 
		We call this quotient the {\bf contracted tensor product} of $M^1$ and $M^2$, and denote the image of $v^1\o v^2$ by $v^1\odot v^2$.
	\end{df}
	The action of $X$ on $M^1\o_\C M^2$ is given by its Laurent expansions at $1$ and $0$, acting on $M^1$ and $M^2$, respectively.
	Using the above spanning elements of $\mathcal{L}_{\mathbb P^1\bs\{1,0\}}(V)_{<0}$, it can be proved that $M^1\odot M^2$ is spanned by symbols $v^1\odot v^2$, bilinear in $v^1$ and $v^2$, subject to the following relations:
	\begin{equation}\label{eq:odot-relation}
		\sum_{j\geq0}\binom{\wt a-s}{j}(a_{(j-t)}v^1)\odot v^2
		=-\sum_{j\geq0}\binom{-t}{j}(-1)^{-t-j}\,v^1\odot(a_{(\wt a-s+j)}v^2),
	\end{equation}
	for all homogeneous $a\in V$, integers $s,t$ with $s+t>1$, and $v^i\in M^i$; see \cite[(4.12)]{Liu26}.
	
	Since the modules in $\Oord$ are generated by their bottom degrees, this quotient is a left $\A=A(L_\ell(\g))$-module by \cite[Theorem~4.8]{Liu26}, with action
	\begin{equation}\label{eq:odot-action}
		[a].(v^1\odot v^2)=\sum_{j\geq0}\binom{\wt a-1}{j}(a_{(j)}v^1)\odot v^2+v^1\odot o(a)v^2,
	\end{equation}
	for homogeneous $a\in V$ and $v^i\in M^i$.
	
	We use the bottom-degree projection $\pi_0$ in \eqref{eq:HL-bottom-projection}. The following lemma is the case $w=1$ of \cite[Theorem~3.14]{CL26}. We sketch the proof for modules in $\Oord$ for the self-containment.
	
	\begin{lm}\label{thm:adm-bottom}
		Let $M^1,M^2\in\Oord$. We have a natural isomorphism of left $\A$-modules
		\begin{equation}\label{eq:etaM1M2}
			\eta_{M^1,M^2}:M^1\odot M^2\xrightarrow{\simeq}\Om(M^1\bt M^2),	\qquad v^1\odot v^2\mapsto\pi_0(v^1\bt v^2),
		\end{equation}
		where $v^1\bt v^2=\mathcal Y_{12}(v^1,e^0)v^2\in\overline{M^1\bt M^2}$, and $\mathcal Y_{12}$ is the canonical intertwining operator in Definition~\ref{df:HL-tensor-product} with $w=1$ and $\log(1)=0$.
	\end{lm}
	\begin{proof}
		Let $W=M^1\bt M^2$. For $M^3\in\Oord$, let $S=\Om(M^3)$. By Lemma~\ref{lm:adm-semisimple}, the contragredient module $(M^3)'$ belongs to $\Oord$ with bottom degree $S^\ast$, which carries the dual right $\A$-action at $\infty$. Since $(M^3)'\cong\Phi^\L(\Om((M^3)'))$ by Lemma~\ref{lm:adm-Zhu-equivalence}, the restriction to $S^\ast$ at $\infty$ is an isomorphism by \cite[Theorem~5.5]{Liu26}.
		
		Denote by $\mathscr C$ and $\mathscr C_{\mathrm{res}}$ the spaces of conformal blocks and $\infty$-restricted conformal blocks on $(\mathbb P^1,\infty,1,0)$, respectively. We have the following commutative diagram of isomorphisms
		\begin{equation}\label{eq:adm-bottom-blocks}
			\begin{tikzcd}[column sep=6em,row sep=2.8em]
				I\fusion{M^1}{M^2}{M^3}
				\arrow[r,"{\text{\cite[Prop.~2.12]{Liu26}}}"] &
				\mathscr C\bigl((M^3)'\o M^1\o M^2\bigr)
				\arrow[d,"\mathrm{res}_{\infty}"',"{\text{\cite[Thm.~5.5]{Liu26}}}"] \\
				\Hom_{\Oord}(W,M^3)
				\arrow[u,"{\text{\cite{HL95a,HL95b}}}","{g\mapsto g\circ\mathcal Y_{12}}"']
				\arrow[d,"\Om","{\text{\cite[Thm.~3.2]{AL26}}}"']
				\arrow[r,dashed] &
				\mathscr C_{\mathrm{res}}\bigl(S^\ast\o M^1\o M^2\bigr)
				\arrow[d,"{\text{\cite[Prop.~4.4, (4.27)]{Liu26}}}"] \\
				\Hom_{\A}(\Om(W),S)
				\arrow[r,dashed,"\beta_{M^3}"] &
				\Hom_{\A}(M^1\odot M^2,S),
			\end{tikzcd}
		\end{equation}
		where, for $g\in 	\Hom_{\Oord}(W,M^3)$, the corresponding restricted block sends $u'\o v^1\o v^2$ to $\langle u'\mid \bar g(v^1\bt v^2)\rangle$ for $u'\in S^\ast$. Then by \eqref{eq:HL-projection-identities}, we have
		\[
		\beta_{M^3}(\Om(g))(v^1\odot v^2)
		=\pi_0\bar g(v^1\bt v^2)
		=\Om(g)\bigl(\pi_0(v^1\bt v^2)\bigr).
		\]
		Taking $M^3=W=M^1\bt M^2$ and $g=\Id_W$, we obtain \eqref{eq:etaM1M2} for $\eta_{M^1,M^2}=\beta_W(\Id_{\Om(W)})$. Thus the bottom arrow sends $f$ to $f\circ\eta_{M^1,M^2}$.
		
		Note that $M^1\odot M^2$ is finite dimensional, and every finite-dimensional left $\A$-module is isomorphic to some $\Om(M^3)$ by Lemma~\ref{lm:adm-Zhu-equivalence}. Hence the desired isomorphism follows from diagram \eqref{eq:adm-bottom-blocks} and the Yoneda lemma. 
	\end{proof}
	
	The following lemma was proved for various classes of affine VOA modules in \cite{Liu26,AL26,CL26}. We briefly sketch its proof. 
	
	\begin{lm}\label{lm:adm-span}
		Let $\la,\mu\in \mathrm{Pr}^{\ell}_{\Z}$. The canonical map
		\begin{equation}\label{eq:varphilamu}
			\varphi_{\la,\mu}:L(\la)\o_\C L(\mu)\ra\Lh{\la}\odot\Lh{\mu},\quad u\o w\mapsto u\odot w,
		\end{equation}
		is a surjective homomorphism of $\g$-modules, where $\g$ acts diagonally on $L(\la)\o_\C L(\mu)$ and through $\A\simeq U(\g)/J_\ell$ on $\Lh{\la}\odot\Lh{\mu}$ via \eqref{eq:odot-action}. 
	\end{lm}
	\begin{proof}
		
		For $x\in\g$, $u\in L(\la)$ and $w\in L(\mu)$, it follows from \eqref{eq:odot-action} that
		\[[x(-1)\vac].(u\odot w)=x(0)u\odot w+u\odot x(0)w.\]
		Hence $\varphi_{\la,\mu}$ is a homomorphism of $\g$-modules.
		
		Note that the modules $\Lh{\la}$ and $\Lh{\mu}$ are spanned by vectors of the form $x^1(-n_1)\cdots x^r(-n_r)u$, where $x^i\in\g$, $n_i\geq1$ and $u$ belongs to $L(\la)$ and $L(\mu)$, respectively. For spanning elements $v^1\in\Lh{\la}$ and $v^2\in\Lh{\mu}$, we use induction on $\deg v^1+\deg v^2$ to show that $v^1\odot v^2\in \mathrm{Im}(\varphi_{\la,\mu})$. The base case when both degrees are zero is clear.
		
		If $\deg v^1>0$, we may assume $v^1=x(-n)w^1$, where $x\in\g$, $n\geq1$ and $\deg w^1=\deg v^1-n$. Then by relation \eqref{eq:odot-relation} with $a=x(-1)\vac$, $s=1$ and $t=n$, we have
		\[x(-n)w^1\odot v^2 =-\sum_{j\geq0}\binom{-n}{j}(-1)^{-n-j}\,w^1\odot x(j)v^2.
		\]
		Note that the sum is finite, and in each nonzero term on the right, the sum of the degrees of the two factors is $\deg w^1+\deg v^2-j<\deg v^1+\deg v^2$.
		If $\deg v^1=0$ and $\deg v^2>0$, we may assume $v^2=x(-n)w^2$ with $n\geq1$. Applying \eqref{eq:odot-relation} with $s=n+1$ and $t=0$, and using $x(j)v^1=0$ for $j\geq1$, we have 
		\[
		v^1\odot x(-n)w^2=-x(0)v^1\odot w^2.
		\]
		Since $x(0)$ preserves degree, the sum of the degrees of the two factors on the right is $\deg v^2-n<\deg v^1+\deg v^2$. In either case, the induction hypothesis proves the claim.
	\end{proof}
	
	For $\la,\mu\in \mathrm{Pr}^{\ell}_{\Z}$, write $K_{\la,\mu}=\ker(\varphi_{\la,\mu})$ and
	\begin{equation}\label{eq:Tlamu}
		T_{\la,\mu}=\frac{L(\la)\o_\C L(\mu)}{K_{\la,\mu}},\qquad 
		\pi_{\la,\mu}:L(\la)\o_\C L(\mu)\twoheadrightarrow T_{\la,\mu}.
	\end{equation}
	Since $\varphi_{\la,\mu}$ in Lemma~\ref{lm:adm-span} is $\g$-linear, it induces an isomorphism of $\g$-modules
	\begin{equation}\label{eq:adm-T}
		\tilde\varphi_{\la,\mu}:T_{\la,\mu}\xrightarrow{\simeq}\Lh{\la}\odot\Lh{\mu},
		\qquad \pi_{\la,\mu}(u\o w)\mapsto u\odot w.
	\end{equation}
	Using this isomorphism, we define the left $\A$-action on $T_{\la,\mu}$ by 
	\[
	[a].\pi_{\la,\mu}(u\o w)=\tilde\varphi_{\la,\mu}^{-1}\big([a].(u\odot w)\big),\qquad [a]\in\A,\quad u\in L(\la),\ w\in L(\mu). 
	\]
	Then $\tilde\varphi_{\la,\mu}$ is an $\A$-module isomorphism. In particular, for $x\in\g$, we have
	\[
	[x(-1)\vac].\pi_{\la,\mu}(u\o w)=\pi_{\la,\mu}(xu\o w+u\o xw).
	\]
	
	\subsection{Explicit braiding formula on the bottom degrees and its eigenvalues}
	
	Using the isomorphism $\eta_{M^1,M^2}$ in Lemma~\ref{thm:adm-bottom}, we can describe the restriction of the braiding \eqref{eq:HL-braiding} to the bottom degrees. 	The following explicit braiding formula with respect to $\odot$-product was proved in  \cite{CL26} in a more general setting. Here we give an alternative proof for modules in the category $\Oord$.
	\begin{prop}\label{thm:adm-contracted-braiding}
		Let $M^1,M^2\in\Oord$. We have an isomorphism of left $\A\simeq U(\g)/J_\ell$-modules
		\begin{equation}\label{eq:R}
			\begin{aligned}
				\mathcal R_{M^1,M^2}:&M^1\odot M^2\xrightarrow{\simeq}M^2\odot M^1,\\
				\mathcal R_{M^1,M^2}(v^1\odot v^2)
				&=e^{\pi i[\om]}.\left(e^{-\pi iL(0)}v^2\odot e^{-\pi iL(0)}v^1\right),
			\end{aligned}
		\end{equation}
		where $[\om]$ acts on $M^2\odot M^1$ by \eqref{eq:odot-action}. Moreover, the following diagram commutes
		\begin{equation}\label{eq:firstbraidingdiagram}
			\begin{tikzcd}[column sep=5em]
				\Om(M^1\bt M^2)\arrow[r,"\Om(c_{M^1,M^2})"] \arrow[d,"\eta_{M^1,M^2}^{-1}"',"\simeq"]& \Om(M^2\bt M^1)\arrow[d,"\eta_{M^2,M^1}^{-1}","\simeq"']\\
				M^1\odot M^2\arrow[r,dashed,"\mathcal R_{M^1,M^2}"]& M^2\odot M^1,
			\end{tikzcd}
		\end{equation}
		i.e., $	\eta_{M^2,M^1}\circ\mathcal R_{M^1,M^2}
		=\Om(c_{M^1,M^2})\circ\eta_{M^1,M^2},$
		where $c_{M^1,M^2}$ is the braiding in the category $\Oord$. 
	\end{prop}
	\begin{proof}
		With the isomorphism $\eta$ in Lemma~\ref{thm:adm-bottom}, we define
		\begin{equation}\label{eq:defofR}
			\mathcal R_{M^1,M^2}
			=\eta_{M^2,M^1}^{-1}\circ\Om(c_{M^1,M^2})\circ\eta_{M^1,M^2}.
		\end{equation}
		Then $\mathcal R_{M^1,M^2}$ is an isomorphism of left $\A$-modules. It remains to show that this map is given by the formula \eqref{eq:R}. By Lemma~\ref{lm:adm-semisimple}, every module in $\Oord$ is a finite direct sum of irreducible modules, it suffices to consider $M^1=\Lh{\la}$ and $M^2=\Lh{\mu}$, with $\la,\mu\in\mathrm{Pr}^{\ell}_{\Z}$.
		Consider the canonical intertwining operator in
		Subsection~\ref{subsec:HL}: 
		\[
		\mathcal Y_{21}\in
		I\fusion{\Lh{\mu}}{\Lh{\la}}{\Lh{\mu}\bt\Lh{\la}},\quad \mathcal Y_{21}(v^2,e^{0})v^1
		=v^2\bt_{P(1)}v^1. 
		\]
		We first determine the bottom-degree component of $\mathcal Y_{21}(v^2,e^{\pi i})v^1$, in view of \eqref{eq:HL-braiding}. 
		
		We may assume $v^1\in\Lh{\la}(r)$ and $v^2\in\Lh{\mu}(s)$, with $L(0)v^1=(h_\la+r)v^1$ and $L(0)v^2=(h_\mu+s)v^2$. By the semisimplicity of $\Oord$, we have a decomposition
		\[
		\Om(\Lh{\mu}\bt\Lh{\la})
		\cong\bigoplus_{\nu\in\mathrm{Pr}^{\ell}_{\Z}}m_\nu L(\nu),
		\qquad m_\nu\in\Z_{\geq0}.
		\] 
		For each $\nu\in\mathrm{Pr}^{\ell}_{\Z}$, let $p_\nu$ denote
		the projection onto the $L(\nu)$-isotypic component, and define
		\[
		\pi_{0,\nu}: \overline{\Lh{\mu}\bt\Lh{\la}}\xrightarrow{\pi_0} \Om(\Lh{\mu}\bt\Lh{\la})=\bigoplus_{\nu\in\mathrm{Pr}^{\ell}_{\Z}}m_\nu L(\nu) \xrightarrow{p_\nu} m_\nu L(\nu). 
		\]
		By the $L(0)$-conjugation formula \cite[(5.4.22)]{FHL}, we have
		\[
		\mathcal Y_{21}(v^2,e^{\pi i})v^1
		=e^{\pi iL(0)}
		\mathcal Y_{21}(e^{-\pi iL(0)}v^2,e^0)
		e^{-\pi iL(0)}v^1.
		\]
		Since $L(0)$ acts on $m_\nu L(\nu)$ by $h_\nu$ and $\pi_{0,\nu}$ commutes with $L(0)$, it follows that
		\[
		\pi_{0,\nu}\left(\mathcal Y_{21}(v^2,e^{\pi i})v^1\right)
		=e^{\pi i(h_\nu-h_\la-h_\mu-r-s)}
		\pi_{0,\nu}\left(\mathcal Y_{21}(v^2,e^0)v^1\right)=e^{\pi i(h_\nu-h_\la-h_\mu-r-s)}
		\pi_{0,\nu}(v^2\bt v^1).
		\]
		Note that $\pi_{0,\nu}e^{L(-1)}=\pi_{0,\nu}$ by \eqref{eq:HL-projection-identities}. Moreover, since $c_{M^1,M^2}$ is a $V$-module homomorphism, we have
		$
		\Om(c_{M^1,M^2})\circ\pi_0
		=\pi_0\circ\overline{c_{M^1,M^2}}.
		$
		Then, by the definition of $\mathcal R_{M^1,M^2}$ in \eqref{eq:defofR}, we have
		\begin{align*}
			p_\nu\eta_{M^2,M^1}
			\bigl(\mathcal R_{M^1,M^2}(v^1\odot v^2)\bigr)
			&=p_\nu\Om(c_{M^1,M^2})\eta_{M^1,M^2}(v^1\odot v^2)=p_\nu\Om(c_{M^1,M^2})\pi_0(v^1\bt v^2)\\
			&=\pi_{0,\nu}\overline{c_{M^1,M^2}}(v^1\bt v^2).
		\end{align*}
		By \eqref{eq:HL-braiding} and the preceding $L(0)$-conjugation formula, we also have
		\begin{align*}
			\pi_{0,\nu}\overline{c_{M^1,M^2}}(v^1\bt v^2)
			&=\pi_{0,\nu}e^{L(-1)}\mathcal Y_{21}(v^2,e^{\pi i})v^1=\pi_{0,\nu}\mathcal Y_{21}(v^2,e^{\pi i})v^1\\
			&=e^{\pi i(h_\nu-h_\la-h_\mu-r-s)}\pi_{0,\nu}(v^2\bt v^1).
		\end{align*}
		On the other hand, since $\eta_{M^2,M^1}$ and $p_\nu$ are $\A$-linear, and $[\om]\in \A$ acts on the $L(\nu)$-isotypic component by $h_\nu\Id$. Therefore,
		\begin{align*}
			p_\nu\eta_{M^2,M^1}\left(e^{\pi i[\om]}.\left(e^{-\pi iL(0)}v^2\odot e^{-\pi iL(0)}v^1\right)\right)&=e^{\pi i h_\nu}e^{-\pi i(h_\mu+s)}e^{-\pi i(h_\la+r)}\pi_{0,\nu}(v^2\bt v^1)\\
			&=e^{\pi i(h_\nu-h_\la-h_\mu-r-s)}\pi_{0,\nu}(v^2\bt v^1).
		\end{align*}
		Since $\nu\in\mathrm{Pr}^{\ell}_{\Z}$ is arbitrary and $\bigoplus_{\nu}p_\nu=\Id$, we have 
		\[
		\eta_{M^2,M^1}
		\bigl(\mathcal R_{M^1,M^2}(v^1\odot v^2)\bigr)=\eta_{M^2,M^1}\left(e^{\pi i[\om]}.\left(e^{-\pi iL(0)}v^2\odot e^{-\pi iL(0)}v^1\right)\right). 
		\]
		Then \eqref{eq:R} follows from the injectivity of $\eta_{M^2,M^1}$ by Lemma~\ref{thm:adm-bottom}. 
	\end{proof}
	
	In order to express the braiding formula \eqref{eq:R} in terms of finite-dimensional $\g$-modules, we first determine the action of $[\om]\in\A$ on $\Lh{\la}\odot\Lh{\mu}$. Let $\{x_a\}$ and $\{x^a\}$ be the dual bases in Subsection~\ref{subsec:adm}, consider the following two elements
	\begin{equation}\label{eq:CandOmega}
		C=\sum_a x_ax^a\in U(\g),\qquad
		\Om=\sum_a x_a\o x^a\in U(\g)\o U(\g),
	\end{equation}
	where $C$ is the Casimir element and $\Om$ is the Casimir tensor. 
	Let $\Delta$ be the standard coproduct on $U(\g)$, we have
	\[
	\Delta(C)=C\o1+1\o C+2\Om.
	\]
	Denote the action of $\Om$ on $L(\la)\o_\C L(\mu)$ by
	$\Om_{\la,\mu}$, and let
	\begin{equation}\label{eq:flip}
		\tau:L(\la)\o_\C L(\mu)\ra L(\mu)\o_\C L(\la),
		\qquad \tau(u\o w)=w\o u,
	\end{equation}
	be the swapping operator. Since
	$\sum_a x_a\o x^a=\sum_a x^a\o x_a$, we have
	\[
	\tau\circ\Om_{\la,\mu}=\Om_{\mu,\la}\circ\tau.
	\]
	
	In view of Lemma~\ref{lm:adm-span}, $\Lh{\la}\odot\Lh{\mu}$ is spanned by $u\odot w$, where  $u\in L(\la)=\Lh{\la}(0)$ and $w\in L(\mu)=\Lh{\mu}(0)$.  The action of $[\om]$ on $\Lh{\la}\odot\Lh{\mu}$ is given by the following formula.
	\begin{lm}\label{lm:adm-omega-action}
		Let $\la,\mu\in\mathrm{Pr}^{\ell}_{\Z}$, $u\in L(\la)$ and $w\in L(\mu)$. Then
		\begin{equation}\label{eq:adm-omega-action}
			[\om].(u\odot w)=(h_\la+h_\mu)u\odot w
			+\frac1\kappa\sum_a x^au\odot x_aw.
		\end{equation}
		In particular, $K_{\la,\mu}=\ker(\varphi_{\la,\mu})$ \eqref{eq:varphilamu} is invariant under $\Om_{\la,\mu}$, and the action of $[\om]$ on $T_{\la,\mu}=(L(\la)\o_\C L(\mu))/K_{\la,\mu}$ is induced by $(h_\la+h_\mu)\Id+\kappa^{-1}\Om_{\la,\mu}$.
	\end{lm}
	\begin{proof}
		Applying \eqref{eq:odot-action} with $a=\om$ and using the
		Segal--Sugawara construction \eqref{eq:sugawara}, together with
		$x^a(j)u=0$ for $j\geq1$, we have
		\begin{align*}
			[\om].(u\odot w)
			&=L(-1)u\odot w+L(0)u\odot w+u\odot L(0)w\\
			&=\frac1\kappa\sum_a x_a(-1)x^a(0)u\odot w
			+(h_\la+h_\mu)u\odot w.
		\end{align*}
		Using relation \eqref{eq:odot-relation} with $a=x_a(-1)\vac$ and $s=t=1$, together with $x_a(j)w=0$ for $j\geq1$, we have
		\[
		x_a(-1)x^a(0)u\odot w=-\sum_{j\geq 0}\binom{-1}{j}(-1)^{-1-j} u\odot x_a(j)w=x^a(0)u\odot x_a(0)w.
		\]
		Substituting this identity into the expression for $[\om].(u\odot w)$ proves \eqref{eq:adm-omega-action}.
		
		By the definition of the $\A$-action on $T_{\la,\mu}$ and
		\eqref{eq:adm-omega-action}, for $u\in L(\la)$ and $w\in L(\mu)$,
		we have
		\begin{align*}
			[\om].\pi_{\la,\mu}(u\o w)
			&=\tilde\varphi_{\la,\mu}^{-1}\bigl([\om].(u\odot w)\bigr)\\
			&=\tilde\varphi_{\la,\mu}^{-1}
			\left((h_\la+h_\mu)u\odot w
			+\frac1\kappa\sum_a x^au\odot x_aw\right)\\
			&=(h_\la+h_\mu)\pi_{\la,\mu}(u\o w)
			+\frac1\kappa\sum_a\pi_{\la,\mu}(x^au\o x_aw).
		\end{align*}
		Since $\sum_a x^a\o x_a=\Om$ \eqref{eq:CandOmega}, for any $\al\in L(\la)\o_\C L(\mu)$, we have 
		\[
		[\om].\pi_{\la,\mu}(\al)
		=\pi_{\la,\mu}\left((h_\la+h_\mu)\al
		+\kappa^{-1}\Om_{\la,\mu}.\al\right).
		\]
		Now, if $\al\in K_{\la,\mu}=\ker(\pi_{\la,\mu})$, this identity gives
		$0=[\om].\pi_{\la,\mu}(\al)
		=\kappa^{-1}\pi_{\la,\mu}(\Om_{\la,\mu}.\al).$
		Hence $\Om_{\la,\mu}.\al\in K_{\la,\mu}$, and so $K_{\la,\mu}$ is
		invariant under $\Om_{\la,\mu}$. Passing to the quotient
		$T_{\la,\mu}=(L(\la)\o_\C L(\mu))/K_{\la,\mu}$, we see that the
		action of $[\om]$ on $T_{\la,\mu}$ is induced by
		$(h_\la+h_\mu)\Id+\kappa^{-1}\Om_{\la,\mu}$.
	\end{proof}
	
	Now, with the identification $\tilde \varphi$ \eqref{eq:adm-T}, we can extend the diagram \eqref{eq:firstbraidingdiagram} to the following commutative diagram
	\begin{equation}\label{eq:firstbraidingdiagram2}
		\begin{tikzcd}[column sep=5em]
			\Om(\Lh{\la}\bt\Lh{\mu})
			\arrow[r,"\Om(c_{\Lh{\la},\Lh{\mu}})"]
			\arrow[d,"\eta_{\la,\mu}^{-1}"',"\simeq"]
			&\Om(\Lh{\mu}\bt\Lh{\la})
			\arrow[d,"\eta_{\mu,\la}^{-1}","\simeq"']
			\\ 
			\Lh{\la}\odot\Lh{\mu}
			\arrow[r,dashed,"\mathcal R_{\Lh{\la},\Lh{\mu}}"] \arrow[d,"\tilde{\varphi}_{\la,\mu}^{-1}"',"\simeq"]
			&	\Lh{\mu}\odot\Lh{\la}\arrow[d,"\tilde{\varphi}_{\mu,\la}^{-1}","\simeq"']\\
			T_{\la,\mu}=\dfrac{L(\la)\o_\C L(\mu)}{K_{\la,\mu}}\arrow[r,dashed,"\mathcal R_{\la,\mu}"] & T_{\mu,\la}=\dfrac{L(\mu)\o_\C L(\la)}{K_{\mu,\la}}.
		\end{tikzcd}
	\end{equation}
	For simplicity, write $\eta_{\la,\mu}=\eta_{\Lh{\la},\Lh{\mu}}$. In other words, under the canonical vertical isomorphisms $\tilde{\varphi}^{-1}\circ\eta^{-1}$, the braiding isomorphism $c_{\Lh{\la},\Lh{\mu}}$ restricts to $\mathcal R: T_{\la,\mu}\to T_{\mu,\la}$.

	\begin{thm}\label{prop:adm-braiding}
		Let $\la,\mu\in\mathrm{Pr}^{\ell}_{\Z}$. The isomorphism $\mathcal R_{\la,\mu}:T_{\la,\mu}\ra T_{\mu,\la}$ in diagram \eqref{eq:firstbraidingdiagram2} is given by the following explicit formula
		\begin{equation}\label{eq:adm-braiding}
			\mathcal R_{\la,\mu}\bigl(\pi_{\la,\mu}(u\o w)\bigr)
			=\pi_{\mu,\la}\left(\tau\exp\left(\frac{\pi i}{\ell+h^\vee}\Om_{\la,\mu}\right)(u\o w)\right),
		\end{equation}
		for $u\in L(\la)$ and $w\in L(\mu)$.
	\end{thm}
	\begin{proof}
		Note that $	\mathcal R_{\la,\mu}
		=\tilde\varphi_{\mu,\la}^{-1}\circ\mathcal R_{\Lh{\la},\Lh{\mu}}\circ\tilde\varphi_{\la,\mu}$, where $\mathcal R_{\Lh{\la},\Lh{\mu}}$ is the isomorphism in Proposition~\ref{thm:adm-contracted-braiding}. 
		Since $(h_\la+h_\mu)\Id$ and $\Om_{\mu,\la}$ commute as operators on $L(\mu)\o_\C L(\la)$, it follows from Lemma~\ref{lm:adm-omega-action} that for $u\in L(\la)$ and $w\in L(\mu)$,
		\begin{equation}\label{eq:exp}
			e^{\pi i[\om]}.\pi_{\mu,\la}(w\o u)
			=e^{\pi i(h_\la+h_\mu)}
			\pi_{\mu,\la}\left(\exp\left(\frac{\pi i}{\kappa}\Om_{\mu,\la}\right)(w\o u)\right).
		\end{equation}
		By \eqref{eq:adm-T}, we have $\tilde\varphi_{\mu,\la}^{-1}(w\odot u)=\pi_{\mu,\la}(w\o u)$. Since $\tilde\varphi_{\mu,\la}^{-1}$ is $\A$-linear, it commutes with the action of $e^{\pi i[\om]}$. 
		Now applying the braiding formula \eqref{eq:R} to $u\odot w$, and using $L(0)u=h_\la u$ and $L(0)w=h_\mu w$, we have
		\begin{align*}
			\mathcal R_{\la,\mu}\bigl(\pi_{\la,\mu}(u\o w)\bigr)
			&=\tilde\varphi_{\mu,\la}^{-1}
			\bigl(\mathcal R_{\Lh{\la},\Lh{\mu}}(u\odot w)\bigr)\\
			&=\tilde\varphi_{\mu,\la}^{-1}
			\left(e^{\pi i[\om]}.\left(e^{-\pi iL(0)}w\odot e^{-\pi iL(0)}u\right)\right)\\
			&=e^{-\pi i(h_\la+h_\mu)}\tilde\varphi_{\mu,\la}^{-1}
			\bigl(e^{\pi i[\om]}.(w\odot u)\bigr)\\
			&=e^{-\pi i(h_\la+h_\mu)}e^{\pi i[\om]}.
			\pi_{\mu,\la}(w\o u).
		\end{align*}
		Substituting \eqref{eq:exp}, using $w\o u=\tau(u\o w)$ and $\tau\circ\Om_{\la,\mu}=\Om_{\mu,\la}\circ\tau$ \eqref{eq:flip}, we have
		\begin{align*}
			\mathcal R_{\la,\mu}\bigl(\pi_{\la,\mu}(u\o w)\bigr)
			&=\pi_{\mu,\la}\left(\exp\left(\frac{\pi i}{\kappa}\Om_{\mu,\la}\right)(w\o u)\right)=\pi_{\mu,\la}\left(\exp\left(\frac{\pi i}{\kappa}\Om_{\mu,\la}\right)\tau(u\o w)\right)\\
			&=\pi_{\mu,\la}\left(\tau\exp\left(\frac{\pi i}{\kappa}\Om_{\la,\mu}\right)(u\o w)\right).
		\end{align*}
		Since $\kappa=\ell+h^\vee$, this proves \eqref{eq:adm-braiding}.
	\end{proof}

	Now we apply Theorem~\ref{prop:adm-braiding} to determine the eigenvalues of the braiding on simple summands of multiplicity one, see  Definition~\ref{df:braiding-eigenvalue}. Let $\la\in\mathrm{Pr}^{\ell}_{\Z}$. Note that
	\[
	L(\la)\o_\C L(\la)
	=S^2L(\la)\oplus\bigwedge\nolimits^2L(\la)
	\]
	is the decomposition into the eigenspaces of the flip $\tau$ in
	\eqref{eq:flip}, with eigenvalues $1$ and $-1$, respectively.
	Both eigenspaces are $U(\g)$-submodules under the diagonal action.
	
	\begin{lm}\label{lm:Tlala}
		The kernel $K_{\la,\la}=\ker(\varphi_{\la,\la})$ \eqref{eq:varphilamu} is invariant under $\tau$. In particular, $\tau$ induces
		an involution on $T_{\la,\la}$ \eqref{eq:Tlamu}, and its eigenspace decomposition is
		\begin{equation}\label{eq:Tlala}
			T_{\la,\la}
			=\pi_{\la,\la}\bigl(S^2L(\la)\bigr)
			\oplus
			\pi_{\la,\la}\bigl(\bigwedge\nolimits^2L(\la)\bigr).
		\end{equation}
		Moreover, both summands are $\A$-submodules of $T_{\la,\la}$.
	\end{lm}
	\begin{proof}
		By Lemma~\ref{lm:adm-omega-action}, $K_{\la,\la}$ is invariant under
		$\Om_{\la,\la}$ and hence under $\exp(-\pi i\Om_{\la,\la}/\kappa)$.
		Let $\al \in K_{\la,\la}$. Then 
		$\b=\exp(-\pi i\Om_{\la,\la}/\kappa)\al\in K_{\la,\la}$ and $\pi_{\la,\la}(\b)=0$ \eqref{eq:Tlamu}. 
		Apply Theorem~\ref{prop:adm-braiding}, we have 
		\[
		0=\mathcal R_{\la,\la}\bigl(\pi_{\la,\la}(\b)\bigr)
		=\pi_{\la,\la}\left(\tau\exp\left(\frac{\pi i}{\kappa}
		\Om_{\la,\la}\right)y\right)
		=\pi_{\la,\la}(\tau \al).
		\]
		Hence $\tau \al\in K_{\la,\la}$, and so $\tau$ induces an involution $\bar{\tau}$ on $T_{\la,\la}$, which is a $\g$-homomorphism. 
		
		Since  $\pi_{\la,\la}$ is surjective and $\bar\tau$ acts by $1$ on $\pi_{\la,\la}(S^2L(\la))$
		and by $-1$ on $\pi_{\la,\la}(\bigwedge\nolimits^2L(\la))$, the quotient $T_{\la,\la}$ decomposes as in \eqref{eq:Tlala}. Since $\bar{\tau}$ is a $\g$-homomorphism, both eigenspaces are $\g$-submodules
		of $T_{\la,\la}$, and hence $\A$-submodules as well, since
		$\A\simeq U(\g)/J_\ell$.
	\end{proof}

	Now write $X=\Lh{\la}$ and $Z=\Lh{\nu}$, for some $\nu\in \Pr^\ell_\Z$. Assume that $Z$ occurs with multiplicity one in $X\bt X\in \Oord$. Recall that $(	\Phi^\L\dashv \Om):\mathsf{Mod}_{\mathrm{fd}}(\A)\rightleftarrows\Oord$ are  equivalence of categories  by Lemma~\ref{lm:adm-Zhu-equivalence}, and 
	$$\tilde{\varphi}_{\la,\la}^{-1}\circ \eta_{\la,\la}^{-1}:\Om(X\bt X)\cong X\odot X\cong  T_{\la,\la}$$
	is a left $\A$-module isomorphism by Lemma~\ref{thm:adm-bottom} and \eqref{eq:adm-T}. It follows that $L(\nu)=\Om(Z)$  occurs with multiplicity one in $T_{\la,\la}$. Then by Lemma~\ref{lm:Tlala}, $L(\nu)$ occurs either in $\pi_{\la,\la}(S^2L(\la))$ or $	\pi_{\la,\la}(\bigwedge^2L(\la))$. Let
	\[
	\epsilon_Z=
	\begin{cases}
		1, & \text{if } L(\nu) \text{ occurs in } \pi_{\la,\la}(S^2L(\la)),\\
		-1, & \text{if } L(\nu) \text{ occurs in }
		\pi_{\la,\la}(\bigwedge\nolimits^2L(\la)).
	\end{cases}
	\]
	
	\begin{coro}\label{coro:adm-braiding-eigenvalues}
		With the notation and assumption as above, the eigenvalue of $c_{X,X}$ on $Z$, in the sense of Definition~\ref{df:braiding-eigenvalue}, is
		\begin{equation}\label{eq:adm-braiding-eigenvalue}
			\lambda_Z=\epsilon_Ze^{\pi i(h_\nu-2h_\la)}
			=\epsilon_Z\exp\left(\frac{\pi i\bigl((\nu\mid\nu+2\rho)-2(\la\mid\la+2\rho)\bigr)}{2(\ell+h^\vee)}\right).
		\end{equation}
	\end{coro}
	\begin{proof}
		Choose an embedding $\iota:Z\hookrightarrow X\bt X$. By
		Definition~\ref{df:braiding-eigenvalue}, $c_{X,X}\circ\iota=\lambda_Z\iota$.
		Applying the functor $\Om$, we see that $\Om(c_{X,X})$ acts by $\lambda_Z$ on the $\A$-submodule $\Om(\iota)(L(\nu))\subseteq\Om(X\bt X)$.
		Let $\mu=\la$ in the diagram \eqref{eq:firstbraidingdiagram2}, have a commutative diagram
		\[
		\begin{tikzcd}[column sep=5em,row sep=2.5em]
			\Om(X\bt X) \arrow[r,"{\Om(c_{X,X})}"]
			& \Om(X\bt X) \\
			T_{\la,\la} \arrow[u,"\simeq"',"\eta_{\la,\la}\circ\tilde\varphi_{\la,\la}"] \arrow[r,"{\mathcal R_{\la,\la}}"]
			& T_{\la,\la}. \arrow[u,"\simeq","\eta_{\la,\la}\circ\tilde\varphi_{\la,\la}"']
		\end{tikzcd}
		\]
		Note that
		$\eta_{\la,\la}\circ\tilde\varphi_{\la,\la}$ carries the unique
		copy of $L(\nu)$ in $T_{\la,\la}$ onto $\Om(\iota)(L(\nu))$.
		Thus it suffices to compute $\mathcal R_{\la,\la}$ on this copy of $L(\nu)\ssq T_{\la,\la}$.
		
		Let $t\in L(\nu)\ssq T_{\la,\la}$. By \eqref{eq:Tlala} and the definition
		of $\epsilon_Z$, we have $\bar\tau (t)=\epsilon_Z \cdot t$. Moreover, by Lemma~\ref{lm:adm-omega-action},  
		$[\om]$ acts as $2h_\la \Id +\kappa^{-1}\bar{\Om}_{\la,\la}$  on $T_{\la,\la}$, where $\bar\Om_{\la,\la}$
		denotes the operator on $T_{\la,\la}$ induced by $\Om_{\la,\la}$. But $[\om]$ also acts as  $h_\nu\Id$ on $L(\nu)$. Hence
		\[
		\kappa^{-1}\bar\Om_{\la,\la}t
		=[\om].t-2h_\la t=(h_\nu-2h_\la)t.
		\]
		Hence, by Theorem~\ref{prop:adm-braiding},
		\[
		\mathcal R_{\la,\la}(t)
		=\bar\tau\exp\left(\frac{\pi i}{\kappa}\bar\Om_{\la,\la}\right)t
		=e^{\pi i(h_\nu-2h_\la)}\bar\tau t
		=\epsilon_Ze^{\pi i(h_\nu-2h_\la)}\cdot t.
		\]
		By commutativity of the diagram, we have $\epsilon_Ze^{\pi i(h_\nu-2h_\la)}=\lambda_Z$.
		Substituting the conformal weight formula \eqref{eq:conformal-weight} gives the second equality
		in \eqref{eq:adm-braiding-eigenvalue}.
	\end{proof}
	
	\begin{remark}
		By \cite[Proposition~6.7]{CHY18}, the braided tensor category $\Oord$ has twist $\theta_M=e^{2\pi iL(0)}$. Applying the balancing identity in \eqref{eq:balancing} to the summand $Z$ of $X\bt X$, we obtain
		\[
		\lambda_Z^2=\theta_Z\theta_X^{-2}=e^{2\pi i(h_\nu-2h_\la)}.
		\]
		Hence the twists determine the eigenvalue of $c_{X,X}$ on $Z$ up to a sign.
		However, to determine the sign, we need the explicit braiding formula of
		Theorem~\ref{prop:adm-braiding} or an equivalent calculation of the
		action of $c_{X,X}$ on the bottom degree $\Om(Z)=L(\nu)$.
	\end{remark}

	
\section{Recognition of quantum group categories of type \texorpdfstring{$A$}{A}}\label{sec:recognition}
	In this section, we prove a recognition theorem for ribbon fusion categories with the fusion rings isomorphic to the fusion ring of affine VOA $L_k(\sl_N)$: such a category is ribbon equivalent to $(\CC(\sl_N,q),c^y)$, compatibly with the labels of the simple objects, if the twist of the vector object and the eigenvalues of its self-braiding agree with those in $(\CC(\sl_N,q),c^y)$; see also Theorem~\ref{thm:recognition}. For $k\geq2$, the proof uses the reconstruction theorem of Kazhdan and Wenzl \cite{KW93} through a comparison result for fusion categories with these fusion rules; see also Proposition~\ref{thm:KW-comparison}. For $k=1$, we use the classification of pointed braided fusion categories.

	\subsection{Fusion rings of \texorpdfstring{$L_k(\sl_N)$}{Lk(slN)}}
	
	Let $\g=\sl_N$ with $N\geq2$. Fix a Cartan subalgebra $\h$ and a set of positive roots $\Delta_+$. Denote the simple roots, simple coroots and fundamental weights by $\al_i$, $\al_i^\vee$ and $\om_i$, respectively, for $1\leq i\leq N-1$, so that $\<\om_i,\al_j^\vee\>=\delta_{ij}$. Let $\theta$ be the highest root and $\rho=\frac12\sum_{\al\in\Delta_+}\al$. Take the normalized invariant bilinear form $(\cdot\mid\cdot)$ such that $(\theta\mid\theta)=2$. We use the convention $\om_2=0$ when $N=2$.
	
	Let $P=\bigoplus_{i=1}^{N-1}\Z\om_i$ be the weight lattice of $\sl_N$ and $Q=\bigoplus_{i=1}^{N-1}\Z\al_i$ the root lattice. Note that $P/Q$ is isomorphic to $\Z/N\Z$ by the isomorphism
	\begin{equation}\label{eq:PQ}
		P/Q\xrightarrow{\sim}\Z/N\Z,\qquad
		\sum_{i=1}^{N-1}a_i\om_i+Q\longmapsto
		\sum_{i=1}^{N-1}i\,a_i\pmod N.
	\end{equation}
	For $\la\in P$, denote by $|\la|$ the image of $[\la]$ under this isomorphism. This is called the {\bf $N$-ality} of $\la$ \cite[Section~1.1]{GG07}.
	
	Recall that a weight $\la\in P$ is {\bf dominant integral} if $\<\la,\al_i^\vee\>\geq0$ for all $i$. The set of dominant integral weights is
	\[
	P_+=\left\{\la=\sum_{i=1}^{N-1}a_i\om_i:a_i\in\Z_{\geq0}\right\}.
	\]
	For $k\in\Z_{\geq0}$, define
	\[
	P^k_+=\{\la\in P_+:(\la\mid\theta)\leq k\}
	=\left\{\sum_{i=1}^{N-1}a_i\om_i:a_i\in\Z_{\geq0},\ \sum_{i=1}^{N-1}a_i\leq k\right\}.
	\]
	For $\la\in P_+$, let $L(\la)$ denote the finite-dimensional irreducible $\g$-module of highest weight $\la$; the Casimir element acts on it by $(\la\mid\la+2\rho)$.

	Let $k\geq1$ and $V=L_k(\sl_N)$ be the associated type $A$ affine VOA.
	The irreducible ordinary $V$-modules are $\{\hat L_k(\la):\la\in P^k_+\}$ \cite{FZ92,DLM97}. Then the fusion ring is
	\begin{equation}\label{eq:affine-fusion-ring}
		\mathfrak F(L_k(\sl_N))=\bigoplus_{\la\in P^k_+}\Z[\hat L_k(\la)],\qquad [\hat L_k(\la)]\cdot[\hat L_k(\mu)]=\sum_{\nu\in P^k_+}N_{\la,\mu}^{\nu}[\hat L_k(\nu)],
	\end{equation}
	where $N_{\la,\mu}^{\nu}$ is the dimension of the space of intertwining operators of type $\binom{\hat L_k(\nu)}{\hat L_k(\la)\;\hat L_k(\mu)}$ \cite{FHL}. 
	
	\begin{definition}\label{df:occur}
		Given  an element
		\[
		\sum_{\la\in P^k_+}m_\la[\hat L_k(\la)]\in 	\mathfrak F(L_k(\sl_N)),
		\qquad m_\la\in\Z_{\geq0},
		\]
		we say that $[\hat L_k(\nu)]$ {\bf occurs} in this element if $m_\nu>0$.
		In particular, $[\hat L_k(\nu)]$ occurs in a product of basis elements if its coefficient in the expansion in the distinguished basis is positive.
	\end{definition}
	
	Let
	\[
	c_{\la,\mu}^{\nu}=\dim\Hom_{\sl_N}(L(\la)\o_\C L(\mu),L(\nu))
	\]
	be the usual tensor-product multiplicities for finite-dimensional $\sl_N$-modules. By the affine fusion rules \cite[Proposition~7.2]{Bea96}, we have
	\begin{equation}\label{eq:fusion-bound}
		0\leq N_{\la,\mu}^{\nu}\leq c_{\la,\mu}^{\nu},\qquad \la,\mu,\nu\in P^k_+.
	\end{equation}
	Moreover,
	\begin{equation}\label{eq:affine-fusion-rule}
		N_{\la,\mu}^{\nu}=c_{\la,\mu}^{\nu}\quad\text{if }(\la+\mu\mid\theta)\leq k+1,
		\qquad \la,\mu,\nu\in P^k_+.
	\end{equation}
	Indeed, if $(\la+\mu\mid\theta)\leq k$, the fusion product has the same decomposition as the ordinary tensor product. If $(\la+\mu\mid\theta)=k+1$, its decomposition is obtained by discarding the constituents $L(\nu)$ with $(\nu\mid\theta)=k+1$; see \cite[Examples~7.3(a,b)]{Bea96}.
	
	Let $\la=\om_i$. Since $(\om_i\mid\theta)=1$ for $1\leq i\leq N-1$ and $(\mu\mid \theta)\leq k$ for $\mu\in P^k_+$, it follows from \eqref{eq:affine-fusion-rule} that
	\begin{equation}\label{eq:fundamental-fusion-rule}
		N_{\om_i,\mu}^{\nu}=c_{\om_i,\mu}^{\nu} \quad\text{for all }
		1\leq i\leq N-1,\quad \mu,\nu\in P^k_+.
	\end{equation}
	Thus, multiplication by $[\hat L_k(\om_i)]$ is obtained from the ordinary tensor product with $L(\om_i)$ by discarding the constituents whose highest weights lie outside $P^k_+$. 
	
	\begin{lm}\label{lm:fundamental-generators}\cite[Lemma~9.5]{Bea96}
		The classes $[\hat L_k(\om_i)]$, $1\leq i\leq N-1$, generate $\mathfrak F(L_k(\sl_N))$ as a ring. 
	\end{lm}
	
	Recall that the fundamental irreducible modules $L(\om_i)$ over $\sl_N$ have the following realizations \cite[Section~15.2]{FH91}:
	\begin{equation}\label{eq:Lomega1wedge}
		L(\om_1)\simeq \C^N=E,\quad L(\om_i)=\bigwedge\nolimits^iE,\qquad 2\leq i\leq N-1. 
	\end{equation} 
	\begin{lm}\label{lm:least-power}
		For the fusion ring $\mathfrak F(L_k(\sl_N))$ in \eqref{eq:affine-fusion-ring}, let
		\begin{equation}\label{eq:blambda}
			b(\la)=\sum_{i=1}^{N-1}ia_i\in \Z_{\geq 0}\quad \mathrm{for}\quad \la=\sum_{i=1}^{N-1}a_i\om_i\in P^k_+. 
		\end{equation}
		Then $b(\la)$ is the least integer $n\geq0$ such that $[\hat L_k(\la)]$ occurs in $[\hat L_k(\om_1)]^n$, see Definition~\ref{df:occur}. 
	\end{lm}
	\begin{proof}
		The assertion is clear for $\la=0$, since $[\hat L_k(\om_1)]^0=[\hat L_k(0)]$ is the unit. Suppose that $\la\neq0$.
		
		For $2\leq i\leq N-1$, in view of \eqref{eq:Lomega1wedge}, the exterior product is an epimorphism of $\sl_N$-modules:
		\[
		E\o_\C\bigwedge\nolimits^{i-1}E\twoheadrightarrow\bigwedge\nolimits^iE,\quad u\o \Lambda\mapsto u\wedge \Lambda.
		\]
		Therefore, $L(\om_i)$ occurs in the regular $\sl_N$-tensor $L(\om_1)\o_\C L(\om_{i-1})$, and so $c^{\om_i}_{\om_1,\om_{i-1}}>0$. 
		Since all the fundamental weights belong to $P^k_+$, \eqref{eq:fundamental-fusion-rule} gives
		\[
		N_{\om_1,\om_{i-1}}^{\om_i}=c_{\om_1,\om_{i-1}}^{\om_i}>0.
		\]
		It follows from this identity and an easy induction on $i\geq 1$ that $[\hat L_k(\om_i)]$ occurs in $[\hat L_k(\om_1)]^i$ for every $1\leq i\leq N-1$. 
		
		For $\la\in P^k_+$, write $\la=\om_{i_1}+\cdots+\om_{i_r}$, where $\om_i$ occurs $a_i$ times and $r=\sum_i a_i=(\la\mid\theta)\leq k$. Let $\mu_0=0$ and $\mu_j=\om_{i_1}+\cdots+\om_{i_j}$ for $1\leq j\leq r$. Then each $\mu_j$ is dominant and satisfies $(\mu_j\mid\theta)=j\leq k$, and so $\mu_j\in P^k_+$. 
		Note that $L(\mu_j)$ occurs exactly once in the $\sl_N$ tensor product module $L(\om_{i_j})\o_\C L(\mu_{j-1})$ as $\om_{i_j}+\mu_{j-1}=\mu_j$. 
		It follows from \eqref{eq:fundamental-fusion-rule} that 
		\begin{equation}\label{eq:Nandc}
			N_{\om_{i_j},\mu_{j-1}}^{\mu_j}=c_{\om_{i_j},\mu_{j-1}}^{\mu_j}=1,
			\qquad 1\leq j\leq r.
		\end{equation}
		If $r=1$, the required occurrence was proved above. Suppose that $r\geq2$. Taking $j=r$ with $\mu_r=\la$ in \eqref{eq:Nandc}, we see that $[\hat L_k(\la)]$ occurs in $[\hat L_k(\om_{i_r})]\cdot [\hat L_k(\mu_{r-1})]$. Then let $j=r-1$ in \eqref{eq:Nandc}, we have $[\hat L_k(\mu_{r-1})]$ occurs in $[\hat L_k(\om_{i_{r-1}})]\cdot [\hat L_k(\mu_{r-2})]$. Hence $[\hat L_k(\la)]$ occurs in $[\hat L_k(\om_{i_r})]\cdot[\hat L_k(\om_{i_{r-1}})]\cdot [\hat L_k(\mu_{r-2})]$. Proceeding inductively, we see that $[\hat L_k(\la)]$ occurs in $[\hat L_k(\om_{i_1})]\cdots [\hat L_k(\om_{i_r})]$. Moreover, each factor $[\hat L_k(\om_{i_j})]$ occurs in $[\hat L_k(\om_1)]^{i_j}$ by the preceding paragraph.  Since all coefficients are nonnegative, $[\hat L_k(\la)]$ therefore occurs in $[\hat L_k(\om_1)]^{i_1+\cdots+i_r}=[\hat L_k(\om_1)]^{b(\la)}$.
		
		Conversely, to show $b(\la)$ is the smallest integer for the occurrence of $[\hat L_k(\la)]$ in $[\hat L_k(\om_1)]^{n}$, we show by induction on $n$ that if $[\hat L_k(\la)]$ occurs in $[\hat L_k(\om_1)]^n$, then $L(\la)$ occurs in $E^{\o n}$, the tensor product of $\sl_N$-modules. The case $n=0$ is clear. For $n\geq1$, write $[\hat L_k(\om_1)]^{n-1}=\sum_{\mu\in P^k_+}m_\mu [\hat L_k(\mu)]$, where $m_\mu\in\Z_{\geq0}$, then 
		\[
		[\hat L_k(\om_1)]^{n}=\sum_{\mu\in P^k_+}m_\mu [\hat L_k(\om_1)]\cdot [\hat L_k(\mu)]= \sum_{\la\in P^k_+}\sum_{\mu\in P^k_+}m_\mu N_{\om_1,\mu}^{\la} [\hat L_k(\la)].
		\]
		
		If $[\hat L_k(\la)]$ occurs in $[\hat L_k(\om_1)]^n$, then $\sum_{\mu\in P^k_+}m_\mu N_{\om_1,\mu}^{\la}>0$, and so $m_\mu N_{\om_1,\mu}^{\la}>0$ for some $\mu\in P^k_+$. Since $m_\mu>0$, $[\hat L_k(\mu)]$ occurs in $[\hat L_k(\om_1)]^{n-1}$. Then by the induction hypothesis, $L(\mu)$ occurs in $E^{\o(n-1)}$, and \eqref{eq:fusion-bound} gives $c_{\om_1,\mu}^{\la}\geq N_{\om_1,\mu}^{\la}>0$. Hence $L(\la)$ occurs in $E\o_\C L(\mu)$ and therefore in $E^{\o n}$.
		
		Choose a basis $e_1,\ldots,e_N$ of $E$ in which $\h$ is diagonal and $\Delta_+$ is the standard positive system, and let $D=\operatorname{diag}(1,\ldots,1,1-N)\in\h$. Then $De_i=e_i$ for $1\leq i\leq N-1$ and $De_N=(1-N)e_N$.  In particular, 
		under the diagonal action of $\sl_N$ on the tensor product module $E^{\o n}$, a standard tensor basis vector $e_{i_1}\o\ds \o e_{i_n}$ with $s$ factors equal to $e_N$ has eigenvalue $(n-s)+s(1-N)=n-Ns\leq n$ under $D$. Hence the eigenvalues of $D$ on $E^{\otimes n}$ are at most $n$. 
		On the other hand, $\om_i(D)=i$ for $1\leq i\leq N-1$, so a highest-weight vector of $L(\la)$ has $D$-eigenvalue $\la(D)=\sum_i ia_i=b(\la)$ \eqref{eq:blambda}. Since $L(\la)$ occurs in $E^{\o n}$, we have $b(\la)\leq n$. 
	\end{proof}
	
	The following proposition will be used to show that an equivalence is compatible with the identifications. It was essentially proved in \cite[Section~4.2]{Gan02}. Here we include an alternative proof for type $A$ using Lemma~\ref{lm:least-power}.
	
	\begin{prop}\label{lm:labels}
		Let $\varphi: \mathfrak F(L_k(\sl_N))\to \mathfrak F(L_k(\sl_N))$ be an automorphism of the fusion ring \eqref{eq:affine-fusion-ring} that permutes the distinguished basis $\{[\hat L_k(\nu)]:\nu \in P^k_+\}$ and fixes the basis element $[\hat L_k(\omega_1)]$. Then $\varphi=\Id$. 
	\end{prop}
	\begin{proof}
		Since $\varphi$ permutes the basis of $ \mathfrak F(L_k(\sl_N))$, there is a permutation $\si:P^k_+\ra P^k_+$ such that
		\[
		\varphi([\hat L_k(\la)])=[\hat L_k(\si(\la))],\qquad \la\in P^k_+.
		\]
		By assumption, $\varphi([\hat L_k(\om_1)])=[\hat L_k(\om_1)]$ and so $\varphi([\hat L_k(\om_1)]^i)=[\hat L_k(\om_1)]^i$ for all $i\geq 0$. Thus $[\hat L_k(\la)]$ occurs in $[\hat L_k(\om_1)]^n$ if and only if $[\hat L_k(\si(\la))]$ occurs in $[\hat L_k(\om_1)]^n$.  Then by Lemma~\ref{lm:least-power}, we have $b(\si(\la))=b(\la)$ for every $\la\in P^k_+$.
		
		Note that the category of ordinary $L_k(\g)$-modules $\mathsf{Ord}(L_k(\g))$ is a rigid braided tensor category, since $L_k(\g)$ is a strongly rational VOA \cite{H08}. The dual object is given by
		\[
		\hat L_k(\la)^\vee\cong\hat L_k(\la^\vee),
		\qquad \la^\vee=-w_0(\la),
		\]
		where $w_0$ is the longest element of the Weyl group such that $w_0(\Delta_+)=-\Delta_+$. Then by the rigidity of the category, we have
		\begin{align*}
			N_{\la,\mu}^{0}
			&=\dim\Hom_{\mathsf{Ord}(L_k(\g))}\left(\hat L_k(\la)\boxtimes\hat L_k(\mu),\vac \right)= \dim\Hom_{\mathsf{Ord}(L_k(\g))}\left(\hat L_k(\mu),\hat L_k(\la)^\vee\right)\\
			&= \dim\Hom_{\mathsf{Ord}(L_k(\g))}\left(\hat L_k(\mu),\hat L_k(\la^\vee)\right).
		\end{align*}
		In particular, in the fusion ring $ \mathfrak F(L_k(\sl_N))$, the unit $[\hat L_k(0)]$ occurs in the product $[\hat L_k(\la)]\cdot[\hat L_k(\mu)]$ if and only if $\mu=\la^\vee$. Since $\varphi([\hat L_k(0)])=[\hat L_k(0)]$, the unit occurs in $\varphi([\hat L_k(\la)]\cdot[\hat L_k(\la^\vee)])=[\hat L_k(\si(\la))]\cdot[\hat L_k(\si(\la^\vee))]$, and so $\si(\la^\vee)=\si(\la)^\vee$. 
		
		For $\la=\sum_{i=1}^{N-1}a_i\om_i$, since $-w_0(\om_i)=\om_{N-i},$ we have $\la^\vee=-w_0(\la)=\sum_{i=1}^{N-1}a_i\om_{N-i}$. It follows from \eqref{eq:blambda} that 
		\[
		\frac{b(\la)+b(\la^\vee)}{N}=
		\frac{\sum_{i=1}^{N-1}ia_i+\sum_{i=1}^{N-1} (N-i)a_i }{N}=\frac{\sum_{i=1}^{N-1}(i+N-i)a_i}{N}
		=\sum_{i=1}^{N-1}a_i=(\la\mid\theta).
		\]
		Now using the fact that $\si$ preserves $b$ and duality, we have 
		\[
		(\si(\la)\mid \theta)=\frac{b(\si(\la))+b(\si(\la)^\vee)}{N}=\frac{b(\la)+b(\la^\vee)}{N}=(\la\mid \theta). 
		\]
		Note that the fundamental dominant weights $\om_i$ are precisely the weights $\la$ in $P^k_+$ such that $(\la\mid \theta)=1$. Hence $\si$ permutes $\{\om_i:1\leq i\leq N-1\}$.  But $b(\om_i)=i=b(\si(\om_i))$ for $1\leq i\leq N-1$ \eqref{eq:blambda}. Hence $\si(\om_i)=\om_i$ and $\varphi([\hat L_k(\om_i)])=[\hat L_k(\om_i)]$ for all $i$. Moreover, since the fusion ring  $\mathfrak F(L_k(\sl_N))$ is generated by $\{[\hat L_k(\om_i)]:1\leq i\leq N-1\}$ by Lemma~\ref{lm:fundamental-generators}, we have $\varphi=\Id$. 
	\end{proof}
	
	\begin{df}
		Let $\mathcal B$ and $\mathcal B'$ be fusion categories whose
		fusion rings are identified, as based rings, with that of
		$L_k(\sl_N)$. Denote by $X_\la$ and $X'_\la$ their
		simple objects corresponding to $\la\in P^k_+$, respectively. 
		A tensor functor $F:\mathcal B\ra\mathcal B'$ is said to be
		{\bf compatible with the identifications} if
		$F(X_\la)\cong X'_\la$ for every $\la\in P^k_+$.
	\end{df}
	
	\subsection{Quantum group categories of type \texorpdfstring{$A$}{A}}\label{sec:QG}
	Throughout this section, we assume that 
	\begin{equation}\label{eq:order-q-square}
		q\in \C\ \text{is a root of unity satisfying}\ 	\mathrm{ord}(q^2)=p\geq N,\quad k=p-N.	
	\end{equation}
	Denote by $U_q(\sl_N)$ Lusztig's quantum group with divided powers at $q$ \cite{Lus93}. We use the comultiplication
	\[
	\Delta(E_i)=E_i\o K_i+1\o E_i,\qquad
	\Delta(F_i)=F_i\o1+K_i^{-1}\o F_i,\qquad
	\Delta(K_i)=K_i\o K_i,
	\]
	where $K_i$ acts on a vector of weight $\la$ by $q^{(\la\mid\al_i)}$.
	
	\begin{definition}[{\cite{A92}}]\label{df:tilting-module}
		Let $M$ be a finite-dimensional type~$1$ $U_q(\sl_N)$-module. A {\bf Weyl filtration} of $M$ is a filtration by submodules
		\[
		0=M_0\subseteq M_1\subseteq\cdots\subseteq M_r=M
		\]
		such that each $M_j/M_{j-1}$ is a Weyl module of highest weight $\la_j\in P_+$. We say that $M$ is a {\bf tilting module} if both $M$ and its dual $M^\ast$ admit a Weyl filtration \cite[Section~3.18]{AP95}.
	\end{definition}
	
	Recall that the tilting modules form a monoidal category $\mathrm{Tilt}(U_q(\sl_N))$  \cite{A92}, with the spherical structure given by $K_{2\rho}$. i.e., for each tilting module $U$, the corresponding pivotal isomorphism is
	\begin{equation}\label{eq:K2rho}
		j_U:U\xrightarrow{\sim}U^{\ast\ast},\qquad
		j_U(u)(\varphi)=\varphi(K_{2\rho}u),
		\qquad u\in U,\quad \varphi\in U^\ast.
	\end{equation}
	Denote by $\mathrm{tr}_q$ the trace for the spherical structure given by $K_{2\rho}$.
	
	\begin{definition}[{\cite{Tur94}}]\label{df:tilting-quotient}
		For tilting modules $U,V$, a morphism $f:U\ra V$ is called  negligible if $\mathrm{tr}_q(g\circ f)=0$ for every morphism $g:V\ra U$. Denote by $\mathcal N(U,V)$ the subspace of negligible morphisms from $U$ to $V$.
		The {\bf semisimplification} of the category $\mathrm{Tilt}(U_q(\sl_N))$ is the quotient
		\begin{equation}\label{eq:tilting-quotient}
			\CC(\sl_N,q):=\mathrm{Tilt}(U_q(\sl_N))/\mathcal N.
		\end{equation}
		The objects of $\CC(\sl_N,q)$ are the tilting modules, and its morphism spaces are
		\[
		\Hom_{\CC(\sl_N,q)}(U,V)=\Hom_{U_q(\sl_N)}(U,V)/\mathcal N(U,V).
		\]
		For a morphism $f:U\ra V$, denote by $[f]$ its class in this quotient.
	\end{definition}
	
	Given tilting modules $U,V\in \operatorname{Tilt}(U_q(\sl_N))$, the tensor product $U\bt V$ in $\CC(\sl_N,q)$ is the image of $U\o_\C V$ with the $U_q(\sl_N)$-action given by the comultiplication $\Delta$.
	Then  $\CC(\sl_N,q)$ \eqref{eq:tilting-quotient} with this tensor product is a fusion category, with simple objects $\{L_q(\nu):\nu\in P^k_+\}$. 
	
	By \cite[Sections~3.20--3.21]{AP95}, for $\la,\mu\in P^k_+$, the tensor products in $\CC(\sl_N,q)$ satisfy
	\begin{equation}\label{eq:QG-fusion-rule}
		L_q(\la)\bt L_q(\mu)\cong\bigoplus_{\nu\in P^k_+}N_{\la,\mu}^{\nu}\,L_q(\nu),
	\end{equation}
	where $N_{\la,\mu}^{\nu}$ are the fusion coefficients of $L_k(\sl_N)$.  In particular, the assignment $[L_q(\la)]\mapsto[\hat L_k(\la)]$ for $\la\in P^k_+$ identifies the fusion ring of $\CC(\sl_N,q)$ with  $\mathfrak F(L_k(\sl_N))$ in \eqref{eq:affine-fusion-ring}.

	\begin{definition}[{\cite{CP94}}]\label{df:universal-R-matrix}
		Let $U,V$ be tilting modules over $U_q(\sl_N)$, and let $y\in\C$ satisfy $y^N=q$. The $R$-matrix on $U\o_\C V$ is the operator
		\begin{equation}\label{eq:QG-universal-R}
			R^y_{U,V}=D^y\prod_{j=s,\dots,1}
			\left(\sum_{m=0}^{p-1}q^{m(m-1)/2}(q-q^{-1})^m[m]_q!\,
			E_{\gamma_j}^{(m)}\o F_{\gamma_j}^{(m)}\right),
		\end{equation}
		where $D^y(u\o w)=y^{N(\al\mid\beta)}u\o w$ for $u\in U_\al$ and $w\in V_\beta$ \cite[(1.3.12)--(1.3.13)]{BK01}, and $\gamma_1,\ldots,\gamma_s$ are the positive roots ordered by a reduced expression for the longest Weyl group element $w_0$, $E_{\gamma_j}^{(m)}$ and $F_{\gamma_j}^{(m)}$ are the corresponding divided powers, and $[m]_q!=\prod_{a=1}^m\frac{q^a-q^{-a}}{q-q^{-1}}$.
	\end{definition}
	For vector spaces $U,V$, let $\tau_{U,V}:U\o_\C V\ra V\o_\C U$ denote the flip, $\tau_{U,V}(u\o w)=w\o u$ as in \eqref{eq:flip}.
	Let $m=\operatorname{ord}(q)$ and choose $y_0\in\C$ of order $Nm$ such that $y_0^N=q$. 
	
	For tilting modules $U,V\in \mathrm{Tilt}(U_q(\sl_N))$, the braiding associated with $y_0$ is $\tau_{U,V}\circ R^{y_0}_{U,V}$ by \cite[Corollary~10.1.20]{CP94}. Passing to the quotient \eqref{eq:tilting-quotient}, we have
	\[
	c^{y_0}_{U,V}=[\tau_{U,V}\circ R^{y_0}_{U,V}]:
	U\bt V\xrightarrow{\simeq}V\bt U,
	\]
	which is the braiding isomorphism on $	\CC(\sl_N,q)$.  
	
	The following lemma describes how the braiding on $\CC(\sl_N,q)$ depend on the choice of an $N$-th root of $q$, which will be used later. 
	A closely related formula appears in \cite[proof of Proposition~2.1.4]{Bru00b}. We include the proof in our notation for completeness.
	
	\begin{lm}\label{lm:QG-braiding}
		Let $\xi\in\C$ satisfy $\xi^N=1$, and let $y=\xi y_0$. Then $y^N=q$. As $U,V$ range over the tilting modules, the following maps
		\begin{equation}\label{eq:QG-braiding}
			c^y_{U,V}=[\tau_{U,V}\circ R^y_{U,V}]:
			U\bt V\xrightarrow{\simeq}V\bt U,
		\end{equation}
		with $R^y_{U,V}$ as in \eqref{eq:QG-universal-R}, define a braiding on $\CC(\sl_N,q)$. If the weights of $U$ and $V$ lie in the classes $a$ and $b$ in $P/Q\cong\Z/N$ \eqref{eq:PQ}, respectively, then
		\begin{equation}\label{eq:QG-braiding-root}
			c^y_{U,V}=\xi^{-ab}c^{y_0}_{U,V}.
		\end{equation}
	\end{lm}
	\begin{proof}
		We prove \eqref{eq:QG-braiding-root} first. Recall that $|\om_i|=i$ under the isomorphism $P/Q\cong\Z/N$. Using the explicit form of the fundamental weights $\om_i$ for $\sl_N$ \cite[Section~13]{Hum}, it is easy to show that $N(\om_i\mid\om_j)=N\min(i,j)-ij$. Hence we have by bilinearity of $(\cdot\mid \cdot)$
		\begin{equation}\label{eq:congruences}
			N(\al\mid\beta)\equiv-|\al|\,|\beta|\pmod N,\qquad \al,\beta\in P.
		\end{equation}
		Then for each pair of weights $\al$ of $U$ and $\beta$ of $V$, we have
		$|\al|=a$ and $|\beta|=b$ by assumption. It follows from \eqref{eq:congruences} that $N(\al\mid \b)\equiv -ab\pmod{N}$, and so
		\[y^{N(\al\mid\beta)}=\xi^{N(\al\mid\beta)}y_0^{N(\al\mid\beta)}
		=\xi^{-ab}y_0^{N(\al\mid\beta)}.
		\]
		Note that only the Cartan factor $D^y$ depends on the choice of $y$ in \eqref{eq:QG-universal-R}. Hence $R^y_{U,V}=\xi^{-ab}R^{y_0}_{U,V}$ on $U\o_\C V$. This shows \eqref{eq:QG-braiding-root} after composing with $\tau_{U,V}$ and passing to the quotient. 
		
		It remains to show $c^y$ satisfies the two hexagon identities \cite{EGNO}. 
		For a tilting module $U$, let $U_\al$ be its weight space of weight $\al$. By grouping the weight spaces according to their classes in $P/Q$, we have
		\[
		U=\bigoplus_{a\in P/Q}U^{(a)},\qquad
		U^{(a)}=\bigoplus_{\al\in P,\ |\al|=a}U_\al.
		\]
		Note that each $U^{(a)}$ is a submodule, since the divided powers $E_i^{(r)}$ and $F_i^{(r)}$ change weights by $r\al_i$ and $-r\al_i$, respectively, and $\al_i\in Q$. By \eqref{eq:QG-braiding-root}, for $a,b\in P/Q$, the map $c^y_{U,V}$ restricts to the isomorphism
		\begin{equation}\label{eq:restrictionofc}
			\xi^{-ab}c^{y_0}_{U^{(a)},V^{(b)}}:
			U^{(a)}\bt V^{(b)}\xrightarrow{\simeq}V^{(b)}\bt U^{(a)}.
		\end{equation}
		
		Let $U,V,W$ be tilting modules, and fix $a,b,c\in P/Q$. Since weights add under the $U_q(\sl_N)$-module tensor product, $V^{(b)}\bt W^{(c)}$ is a summand of $(V\bt W)^{(b+c)}$, and $U^{(a)}\bt V^{(b)}$ is a summand of $(U\bt V)^{(a+b)}$. Since the associator on the fusion category $\CC(\sl_N,q)$ is canonical, by \eqref{eq:restrictionofc} and the first hexagon identity for $c^{y_0}$, we have
		\begin{align*}
			c^y_{U^{(a)},V^{(b)}\bt W^{(c)}}
			&=\xi^{-a(b+c)}c^{y_0}_{U^{(a)},V^{(b)}\bt W^{(c)}}=\xi^{-ac}\xi^{-ab}(\Id_{V^{(b)}}\bt c^{y_0}_{U^{(a)},W^{(c)}})\circ(c^{y_0}_{U^{(a)},V^{(b)}}\bt\Id_{W^{(c)}})\\
			&=(\Id_{V^{(b)}}\bt c^y_{U^{(a)},W^{(c)}})\circ(c^y_{U^{(a)},V^{(b)}}\bt\Id_{W^{(c)}}).
		\end{align*}
		The other hexagon identity for $c^{y}$ is similar, we omit the details.
	\end{proof}
	
	\begin{definition}\label{df:QG-category}
		Given $y\in\C$ with $y^N=q$, let $$(\CC(\sl_N,q),c^y)$$ be the ribbon fusion category $\CC(\sl_N,q)$ with braiding $c^y$ \eqref{eq:QG-braiding} and spherical structure given by $K_{2\rho}$ \eqref{eq:K2rho}. 
	\end{definition}
	\begin{lm}
		The twist on $L_q(\la)$ in the ribbon category $(\CC(\sl_N,q),c^y)$ is given by
		\begin{equation}\label{eq:QG-twist}
			\theta_{L_q(\la)}=y^{N(\la\mid\la+2\rho)} \Id_{L_q(\la)},\qquad \la\in P^k_+.
		\end{equation}
	\end{lm}
	\begin{proof}
		Note that for the usual braiding $c^{y_0}$ on  $	\CC(\sl_N,q)$, the twist on $L_q(\la)$ is $y_0^{N(\la\mid\la+2\rho)}\Id_{L_q(\la)}$ \cite{BK01,Saw06}. 	By \eqref{eq:QG-braiding-root} and the partial trace formula for the twist \cite{BK01,EGNO}, replacing $y_0$ by $y=\xi y_0$ multiplies the twist on $L_q(\la)$ by $\xi^{-|\la|^2}$.
		Since $2\rho\in Q$, we have $\xi^{-|\la|^2}=\xi^{N(\la\mid\la+2\rho)}$ by \eqref{eq:congruences}.
		Therefore, for the braiding $c^y$, we have 
		\[
		\theta_{L_q(\la)}
		=\xi^{-|\la|^2}y_0^{N(\la\mid\la+2\rho)}\Id_{L_q(\la)}
		=(\xi y_0)^{N(\la\mid\la+2\rho)}\Id_{L_q(\la)},
		\]
		which gives \eqref{eq:QG-twist} since $y=\xi y_0$.
	\end{proof}
	
	Let $E=L(\om_1)=\C^N$ as in \eqref{eq:Lomega1wedge}. It follows from \cite[Proposition~15.15]{FH91} that
	\begin{equation}\label{eq:EtensorE}
		E\o_\C E=S^2E\oplus\bigwedge\nolimits^2 E,\qquad
		S^2E\cong L(2\om_1),\quad \bigwedge\nolimits^2 E\cong L(\om_2),
	\end{equation}
	where $\om_2=0$ when $N=2$. Using $(\om_i\mid 2\rho)=i(N-i)$ and $(\om_i\mid\om_j)=\min(i,j)-ij/N$, we obtain
	\begin{equation}\label{eq:casimirs}
		\begin{aligned}
			(\om_1\mid\om_1+2\rho)&=\frac{N^2-1}{N},&
			(2\om_1\mid2\om_1+2\rho)&=\frac{2(N-1)(N+2)}{N},\\
			(\om_2\mid\om_2+2\rho)&=\frac{2(N-2)(N+1)}{N}.
		\end{aligned}
	\end{equation}
	
	For $k\geq1$, let $Y=L_q(\om_1)$. It follows from \eqref{eq:QG-twist} and \eqref{eq:casimirs} that the twist on $Y$ is given by 
	\[
	\theta_Y=y^{N(\om_1\mid\om_1+2\rho)}\Id_Y=y^{N^2-1}\Id_Y.
	\]
	\begin{prop}\label{prop:QG-eigenvalues}
		Let $y^N=q$, and let $c^y$ be the braiding in \eqref{eq:QG-braiding}. With the braiding eigenvalues below are defined in Definition~\ref{df:braiding-eigenvalue}, we have the following
		\begin{enumerate}[label=(\arabic*)]
			\item If $k=0$, then $\CC(\sl_N,q)\simeq\mathrm{Vec}$.
			\item If $k=1$, then $Y\bt Y\cong L_q(\om_2)$ and the eigenvalue of $c^y_{Y,Y}$ on $L_q(\om_2)$ is $\lambda_-=-y^{-N-1}$.
			\item If $k\geq2$, then $Y\bt Y\cong L_q(2\om_1)\oplus L_q(\om_2)$. The respective eigenvalues of $c^y_{Y,Y}$ on $L_q(2\om_1)$ and $L_q(\om_2)$ are given by 
			\begin{equation}\label{eq:QG-eigenvalues}
				\lambda_+=qy^{-1}=y^{N-1},\qquad \lambda_-=-q^{-1}y^{-1}=-y^{-N-1}.
			\end{equation}
		\end{enumerate}
	\end{prop}
	\begin{proof}
		If $k=0$, then $P^0_+=\{0\}$, so the tensor unit is the only simple object.
		
		For $k\geq1$, apply \eqref{eq:QG-fusion-rule} to $\la=\mu=\om_1$. By \eqref{eq:fundamental-fusion-rule}, for each $\nu\in P^k_+$, the multiplicity of $L_q(\nu)$ in $Y\bt Y$ is $N_{\om_1,\om_1}^{\nu}=c_{\om_1,\om_1}^{\nu}$. Since $E\o_\C E\cong L(2\om_1)\oplus L(\om_2)$ \eqref{eq:EtensorE}, the only possible summands are $L_q(\om_2)$ and $L_q(2\om_1)$, each with multiplicity one whenever its label belongs to $P^k_+$. 
		Note that  we have $\om_2\in P^k_+$ for all $k\geq1$, whereas $2\om_1\in P^k_+$ if and only if $k\geq2$.
		
		Then the braiding eigenvalues follow from \cite[Section~3(B)]{Ker94} with $t=y^{-1}$.
	\end{proof}
	
	\subsection{The Kazhdan--Wenzl invariants}\label{subsec:KW-invariants}
	Throughout this subsection, we assume $N,k\geq2$ are integers and $p=N+k$.
	
	We first recall some basic notions in Hecke algebras. We omit the associativity and unit isomorphisms in the formulas below. 
	
	\begin{definition}\label{df:KW-Hecke-parameter}
		Let $\mathcal A$ be a fusion category with fusion ring $\mathfrak F(\mathcal A)$ isomorphic to $\mathfrak F(L_k(\sl_N))$ \eqref{eq:affine-fusion-ring} as a based ring. In the notation of \cite[Section~1]{KW93}, $\mathcal A$ is an $SL_{N,p}$-category. Using such an isomorphism, denote the simple objects of $\mathcal A$ by $Z_\la$, $\la\in P^k_+$, and let $Z=Z_{\om_1}$. By \eqref{eq:fundamental-fusion-rule}, we have
		\[
		Z\bt Z\cong Z_{2\om_1}\oplus Z_{\om_2}.
		\]
		Let $a\in\End_{\mathcal A}(Z\bt Z)$ be the projection onto $Z_{\om_2}$. In particular, for $N=2$, $a$ is the projection onto the unit summand.
		
		For $m\geq2$ and $Q\in\C^\times$, let $H_m(Q)$ be the Hecke algebra with generators $g_1,\ldots,g_{m-1}$ and relations
		\begin{equation}\label{eq:KW-Hecke-relations}
			\begin{aligned}
				g_ig_{i+1}g_i&=g_{i+1}g_ig_{i+1} &&(1\leq i\leq m-2),\\
				g_ig_j&=g_jg_i &&(|i-j|\geq2),\\
				(g_i-Q)(g_i+1)&=0 &&(1\leq i<m).
			\end{aligned}
		\end{equation}
		Define the following operators
		\begin{equation}\label{eq:KW-generator}
			T=Qa-(\Id-a),\qquad
			T_i=\Id_Z^{\bt(i-1)}\bt T\bt\Id_Z^{\bt(m-i-1)},\qquad 1\leq i<m.
		\end{equation}
		Then $T$ acts by $Q\Id $ on $Z_{\om_2}$ and by $-\Id$ on $Z_{2\om_1}$. 
		
		We call a primitive $p$-th root $Q\in \C^\times$ a {\bf Hecke parameter} of $\mathcal A$ if the associated  $T_i$ satisfy \eqref{eq:KW-Hecke-relations} for every $m\geq2$.
	\end{definition}
	
	Note that a Hecke parameter of $\mathcal A$ is not unique. if $Q$ is a Hecke parameter of $\mathcal A$, then so is $Q^{-1}$; see \cite[Proposition~4.2]{KW93}.
	
	With $\mathcal A$ and $Z$ as in Definition~\ref{df:KW-Hecke-parameter}, we define the second invariant. Note that the unit $\vac$ occurs with multiplicity one in $Z^{\bt N}$. Indeed, since $\mathfrak F(\mathcal A)\cong \mathfrak F(L_k(\sl_N))$ by assumption, it follows from Lemma~\ref{lm:least-power} that $Z^*\cong Z_{\om_{N-1}}$ occurs in $Z^{\bt(N-1)}$. Hence $\vac$ occurs in $Z^{\bt N}$. On the other hand, it follows from \eqref{eq:fusion-bound} that the multiplicity of $\vac$ in $Z^{\bt N}$ is at most that of the trivial $\sl_N$-module in $(\C^N)^{\o N}$, which is one. Hence $\vac$ occurs with multiplicity exactly one in $Z^{\bt N}$. 
	\begin{definition}
		Choose morphisms
		\[
		\iota:\vac\ra Z^{\bt N},\qquad \pi:Z^{\bt N}\ra\vac,\qquad \pi\circ\iota=\Id_\vac.
		\]
		For a Hecke parameter $Q$, the {\bf determinant scalar} $\tau_{\mathcal A}$ is given by
		\begin{equation}\label{eq:KW-determinant}
			\tau_{\mathcal A}\Id_Z=(\pi\bt\Id_Z)\circ T_N\cdots T_1\circ(\Id_Z\bt\iota),
		\end{equation}
		where $T_1,\ldots,T_N$ act on $Z^{\bt(N+1)}$ via \eqref{eq:KW-generator}. For the normalization of the determinant scalar in \cite[Proposition~5.1]{KW93}, see \cite[Section~7]{PR11}. 
	\end{definition}
	
	For fixed $Q$, the scalar $\tau_{\mathcal A}$ in \eqref{eq:KW-determinant} is independent of the choice of $\iota$ and $\pi$ \cite[p.~135]{KW93}.

	For a primitive $p$-th root of unity $Q$, let $e_N^+\in H_N(Q)$ be the primitive central idempotent corresponding to the one-dimensional representation on which each generator $g_i$ acts by $Q \Id$.
	
	\begin{lm}[{\cite[Theorem~4.1 and Lemma~4.7]{KW93}}]\label{lm:KW-Hecke}
		With the assumptions and notation of Definition~\ref{df:KW-Hecke-parameter}, $\mathcal A$ admits a Hecke parameter $Q$. For every such $Q$ and every $m\geq2$, the representation morphism
		\[
		\rho_m:H_m(Q)\ra\End_{\mathcal A}(Z^{\bt m}),\qquad g_i\mapsto T_i,
		\]
		is surjective, and its kernel depends only on $N$, $p$ and $Q$. The induced isomorphisms of the quotients are compatible with tensor products. For $m=N$, we have $\rho_N(e_N^+)=\iota\circ\pi$.
	\end{lm}
	In the lemma above, the parameters of Kazhdan and Wenzl are $\ell=p$ and $q_{\mathrm{KW}}=Q^{-1}$ \cite{KW93}. Their Hecke generator is written as $T_{\mathrm{KW}}=-Q^{-1}T$, so their sign character corresponds to $g_i\mapsto Q$ in our notation.
	
	We also recall some definitions and lemmas from \cite{TW05} for arbitrary fusion categories. We will use them in the proof of Proposition~\ref{thm:KW-comparison}.
	\begin{lm}[{\cite[Theorem~3.5]{TW05}}]\label{lm:tensor-extension}
		Let $\mathcal D$ and $\mathcal D'$ be full monoidal subcategories of fusion categories $\mathcal C$ and $\mathcal C'$, respectively. Suppose that every simple object of $\mathcal C$ (resp. $\mathcal C'$) is a direct summand of an object of $\mathcal D$ (resp. $\mathcal D'$). Then every $\C$-linear tensor equivalence of subcategories $\mathcal D\xrightarrow{\simeq}\mathcal D'$ extends to a tensor equivalence of the  whole categories $\mathcal C\xrightarrow{\simeq}\mathcal C'$.
	\end{lm}
	
	\begin{definition}[{\cite[Definitions~4.1, 4.3 and (4.1)]{TW05}}]\label{df:TW-diagonal}
		Let $\mathcal C$ be a fusion category, let $X$ be an object of $\mathcal C$, and let $N\geq1$ be an integer.
		
		Let $\mathcal D$ be the full subcategory of $\mathcal C$ whose objects are $X^{\bt m}$ for $m\geq 0$, where $X^{\bt 0}=\vac$.
		We call $\mathcal D$ the {\bf monoidal algebra generated by $X$}.
		It is said to be {\bf of type $N$} if $X$ is simple, $\vac$ occurs with
		multiplicity one in $X^{\bt N}$, and
		\[
		\Hom_{\mathcal C}(X^{\bt m},X^{\bt n})=0
		\quad\text{whenever }m\not\equiv n\pmod N.
		\]
		The {\bf diagonal} $\Delta\mathcal D$ is a monoidal category that has the same objects as $\mathcal D$, with morphism spaces
		\[
		\Hom_{\Delta\mathcal D}(X^{\bt m},X^{\bt n})=
		\begin{cases}
			\End_{\mathcal C}(X^{\bt m}),&m=n,\\
			0,&m\ne n.
		\end{cases}
		\]
		Assume that $\mathcal D$ is of type $N$ and the category $\Delta\mathcal D$ has a braiding $c$. Choose $\iota:\vac\ra X^{\bt N}$ and $\pi:X^{\bt N}\ra\vac$ with $\pi\circ\iota=\Id_\vac$. The scalar $\Theta_{\mathcal D}$ associated to $c$ is given by
		\begin{equation}\label{eq:TW-diagonal-scalar}
			\Theta_{\mathcal D}\Id_X=(\pi\bt\Id_X)\circ c_{X,X^{\bt N}}\circ(\Id_X\bt\iota).
		\end{equation}
		Two monoidal algebras $\mathcal D$ and $\mathcal D'$ of type $N$, generated by $X$ and $X'$, respectively, are called {\bf extensions of the same braided diagonal} if there is a braided tensor equivalence $\Psi:\Delta\mathcal D\xrightarrow{\simeq}\Delta\mathcal D'$ such that $\Psi(X^{\bt m})=(X')^{\bt m}$ for every $m\geq0$.
	\end{definition}
	
	\begin{lm}[{\cite[Proposition~4.6]{TW05}}]\label{lm:TW-diagonal-extension}
		Let $\mathcal D$ and $\mathcal D'$ be extensions of the same braided diagonal. Let $\iota$ and $\pi$ be the inclusion and projection of the unit summand of $X^{\bt N}$, and let $\iota'$ and $\pi'$ be the corresponding maps for $(X')^{\bt N}$. 
		
		With $\Psi:\Delta\mathcal D\xrightarrow{\simeq}\Delta\mathcal D'$ as in Definition~\ref{df:TW-diagonal}, assume that
		\[
		\Psi(\iota\circ\pi)=\iota'\circ\pi'
		\qquad\text{and}\qquad
		\Theta_{\mathcal D}=\Theta_{\mathcal D'}.
		\]
		Then $\Psi$ extends to a tensor equivalence $\mathcal D\xrightarrow{\simeq}\mathcal D'$ such that $\Psi(X)=X'$.
	\end{lm}
	
	We now return to the categories whose fusion rings are $L_k(\sl_N)$ and give the following form of the Kazhdan--Wenzl reconstruction theorem \cite{KW93}. 
	
	\begin{prop}\label{thm:KW-comparison}
		Let $\mathcal A$ and $\mathcal A'$ be fusion categories with fusion rings $\mathfrak F(\mathcal A)$ and $\mathfrak F(\mathcal A')$ both isomorphic to $\mathfrak F(L_k(\sl_N))$ as based rings. Denote their simple objects by $Z_\la$ and $Z'_\la$, respectively, for $\la\in P^k_+$, and write  $Z=Z_{\om_1}$ and $Z'=Z'_{\om_1}$. 
		
		Let $Q,Q'$ be Hecke parameters of $\mathcal A,\mathcal A'$, respectively, and let $\tau_{\mathcal A},\tau_{\mathcal A'}$ be the corresponding determinant scalars in \eqref{eq:KW-determinant}. 
		If $Q=Q'$ and $\tau_{\mathcal A}=\tau_{\mathcal A'}$, then there is a tensor equivalence $F:\mathcal A\xrightarrow{\simeq}\mathcal A'$ such that $F(Z_\la)\cong Z'_\la$ for every $\la\in P^k_+$.
	\end{prop}
	\begin{proof}
		Let $\mathcal D$ and $\mathcal {D}'$ be the monoidal algebras generated by $Z$ in $\mathcal A$ and by $Z'$ in $\mathcal A'$, respectively. By \eqref{eq:fusion-bound}, every simple summand of $Z^{\bt m}$ has $N$-ality $m\bmod N$. Hence $$\Hom_{\mathcal A}(Z^{\bt m},Z^{\bt n})=0\quad \mathrm{unless}\quad m\equiv n\pmod N,$$ 
		and the same holds for $Z'$. Since the unit occurs with multiplicity one in each $N$-th tensor power, $\mathcal D$ and $\mathcal D'$ are of type $N$ by Definition~\ref{df:TW-diagonal}.
		
		For $m\geq2$, denote by $\rho_m$ and $\rho'_m$ the representations of $H_m(Q)$ on $Z^{\bt m}$ and $(Z')^{\bt m}$, respectively. Since $Q=Q'$, it follows from Lemma~\ref{lm:KW-Hecke} that $\ker\rho_m=\ker\rho'_m$. Thus $\rho_m$ and $\rho'_m$ induce the algebra isomorphisms
		\begin{equation}\label{eq:twoEnd}
			\End_{\mathcal A}(Z^{\bt m})
			\xleftarrow[\sim]{\bar\rho_m} H_m(Q)/\ker\rho_m
			\xrightarrow[\sim]{\bar\rho'_m} \End_{\mathcal A'}((Z')^{\bt m}),
		\end{equation}
		which preserve tensor products. 
		
		\noindent {\bf Step 1.} We construct braidings on diagonals categories $\Delta\mathcal D$ and $\Delta\mathcal D'$, and a braided tensor equivalence $\Psi: \Delta\mathcal D\xrightarrow{\simeq}\Delta\mathcal D'$.

		Choose $q,y\in \C^\times$ such that $q^2=Q^{-1}$ and $y^N=q$. 
		Consider the fusion category $\mathcal A''=\CC(\sl_N,q)$, and let $Y=L_q(\om_1)$.
		Since $k\geq2$, we have $Y\bt Y\cong L_q(2\om_1)\oplus L_q(\om_2)$ by Proposition~\ref{prop:QG-eigenvalues}. Consider the operator in \eqref{eq:KW-generator} \[T=Qa-(\Id-a): Y\bt Y\to Y\bt Y,\]
		where $a: Y\bt Y\to L_q(\om_2)$ is the projection. Then  $T(u)=-u$ for any $u\in L_q(2\om_1)$ and $T(v)=Qv$ for any $v\in L_q(\om_2)$. We claim that $c^y_{Y,Y}=sT$ on $Y\bt Y$, where $s=-y^{N-1}$. 
		
		Indeed, for $u$ and $v$ as above, it follows from \eqref{eq:QG-eigenvalues} and $Q=q^{-2}=y^{-2N}$ that
		\[
		\begin{alignedat}{4}
			sT(u) &= -su=  y^{N-1}u  = c^y_{Y,Y}(u), &\qquad &u\in L_q(2\om_1),\\
			sT(v) &= sQv = -y^{-N-1}v = c^y_{Y,Y}(v), &       &v\in L_q(\om_2).
		\end{alignedat}
		\]
		Thus $c^y_{Y,Y}=sT$ on $Y\bt Y$. Therefore, on $Y^{\bt m}$, each $T_i$ is $s^{-1}$ times the braiding on the $i$-th and $(i+1)$-st factors. It follows that $T_i$ satisfy the first two relations in \eqref{eq:KW-Hecke-relations}. The quadratic relation follows from $(T-Q\Id)(T+\Id)=0$ on $Y\bt Y$.
		Hence $Q$ is also a Hecke parameter of $\mathcal A''$. 
		
		Apply Lemma~\ref{lm:KW-Hecke}  to $\mathcal A''$ with Hecke parameter $Q$, we have $	\End_{\mathcal A''}(Y^{\bt m})\cong H_m(Q)/\ker\rho_m$. Then it follows from \eqref{eq:twoEnd} that there exists natural isomorphism
		\[\al_r	: \End_{\mathcal A''}(Y^{\bt r})\xrightarrow{\simeq}\End_{\mathcal A}(Z^{\bt r}),
		\qquad r\geq2,\]
		which preserves tensor product. 
		Let $\al_0=\al_1=\Id_\C$ and define
		\[
		c_{m,n}:=\al_{m+n}\bigl(c^y_{Y^{\bt m},Y^{\bt n}}\bigr)\in \End_{\mathcal A}(Z^{\bt (m+n)}),\qquad m,n\geq0.
		\]
		Since the $\al_r$ preserve tensor products, applying  $\al_{m+n+r}$ to the hexagon identities for $c^y$ shows the corresponding identities for the $c_{m,n}$. Hence the $c_{m,n}$ define a braiding on $\Delta\mathcal D$. Similarly, let $\al'_r:\End_{\mathcal A''}(Y^{\bt r})\xrightarrow{\simeq}\End_{\mathcal A'}((Z')^{\bt r})$ be the corresponding isomorphism for $\mathcal A'$. Then the maps
		\[
		c'_{m,n}:=\al'_{m+n}\bigl(c^y_{Y^{\bt m},Y^{\bt n}}\bigr)\in  \End_{\mathcal A'}(Z^{\bt (m+n)}),\qquad m,n\geq0,
		\]
		define a braiding on $\Delta\mathcal D'$.
		
		For $r\geq2$, the isomorphism in \eqref{eq:twoEnd} is $\al'_r\circ\al_r^{-1}$. Since these maps preserve tensor products by Lemma~\ref{lm:KW-Hecke},  and are the identity on $\C$ for $r=0,1$, it follows from Definition~\ref{df:TW-diagonal} that we have a tensor equivalence
		\[
		\Psi:\Delta\mathcal D\xrightarrow{\simeq}\Delta\mathcal D',\qquad Z^{\bt r}\mapsto (Z')^{\bt r},\qquad f\mapsto\al'_r\bigl(\al_r^{-1}(f)\bigr).
		\]
		In particular, by the definition of $c_{m,n}$ and $c'_{m,n}$, we have $\Psi(c_{m,n})=\al'_{m+n}\bigl(c^y_{Y^{\bt m},Y^{\bt n}}\bigr)=c'_{m,n}$ for all $m,n\geq0$. Hence $\Psi$ is also a braided equivalence. 
		
		\noindent {\bf Step 2.} We lift $	\Psi:\Delta\mathcal D\xrightarrow{\simeq}\Delta\mathcal D'$ to a tensor equivalence $F: \mathcal A\xrightarrow{\simeq} \mathcal A'$. 
		
		By \eqref{eq:twoEnd}, $\Psi$ acts on $\End_{\mathcal A}(Z^{\bt N})$ by $\bar\rho'_N\circ\bar\rho_N^{-1}$. Hence, by Lemma~\ref{lm:KW-Hecke}, we have
		\[
		\Psi(\iota\circ\pi)
		=\Psi\bigl(\rho_N(e_N^+)\bigr)
		=\rho'_N(e_N^+)
		=\iota'\circ\pi'.
		\]
		By the hexagon identity and $c^y_{Y,Y}=sT$, we have $c_{1,N}=c_{Z,Z^{\bt N}}=s^NT_N\cdots T_1$, and similarly for $c'_{1,N}$. Substituting this braiding formula into \eqref{eq:TW-diagonal-scalar} and using \eqref{eq:KW-determinant}, we have
		\begin{align*}
			\Theta_{\mathcal D}\mathrm{Id}_Z
			&=(\pi\boxtimes\mathrm{Id}_Z)\circ c_{1,N}\circ(\mathrm{Id}_Z\boxtimes\iota)\\
			&=s^N(\pi\boxtimes\mathrm{Id}_Z)\circ T_N\cdots T_1\circ(\mathrm{Id}_Z\boxtimes\iota)\\
			&=s^N\tau_{\mathcal A}\,\mathrm{Id}_Z,
		\end{align*}
		Similarly, $\Theta_{\mathcal D'}\mathrm{Id}_{Z'}=s^N\tau_{\mathcal A' }\,\mathrm{Id}_{Z'}$. Since the determinant scalars $\tau_{\mathcal A}=\tau_{\mathcal A'}$ by assumption, we have 
		\[
		\Theta_{\mathcal D}=s^N\tau_{\mathcal A}=s^N\tau_{\mathcal A'}=\Theta_{\mathcal D'}.
		\]
		Hence $\Psi: \Delta\mathcal D\xrightarrow{\simeq}\Delta\mathcal D'$ satisfies both assumptions in  Lemma~\ref{lm:TW-diagonal-extension}, and so $\Psi$ extends to a tensor equivalence $\mathcal D\xrightarrow{\simeq}\mathcal D'$. Furthermore, every simple object of $\mathcal A$ occurs in some $Z^{\bt m}$ by Lemma~\ref{lm:least-power}, and similarly for $\mathcal A'$ and $Z'$. It follows from Lemma~\ref{lm:tensor-extension} that $\Psi$ extends to a tensor equivalence $F:\mathcal A\xrightarrow{\simeq}\mathcal A'$. 
		
		Finally, since $F(Z_{\om_1})\cong Z'_{\om_1}$, it follows from Proposition~\ref{lm:labels} that the induced map of $F$ on the fusion rings of $\mathcal A$ and $\mathcal A'$ is identity. Hence we have  $F(Z_\la)\cong Z'_\la$ for all $\la\in P^k_+$.
	\end{proof}
	
	\begin{remark}\label{rem:KW93}
		By the reconstruction theorem of Kazhdan and Wenzl \cite[Theorem~$\mathrm{A}_\ell$]{KW93}, $\mathcal A$ is monoidally equivalent to $\CC(\sl_N,q)$ with its associator twisted by a $3$-cocycle on $\Z/N$, for some $q$ with $\operatorname{ord}(q^2)=p$. To deduce Proposition~\ref{thm:KW-comparison} from this statement, one still needs to compare the cocycle classes for $\mathcal A$ and $\mathcal A'$.  We therefore give a direct proof of Proposition~\ref{thm:KW-comparison} using the extension result of Tuba and Wenzl \cite[Proposition~4.6]{TW05}.
	\end{remark}
	
	\subsection{The recognition theorem for type \texorpdfstring{$A$}{A}}\label{subsec:qg-recognition}
	Now we prove the main theorem of this section, which gives sufficient conditions for a ribbon fusion category $\mathcal B$ to be ribbon equivalent to $(\CC(\sl_N,q),c^y)$ in Definition~\ref{df:QG-category}. We first give a lemma for the case $k=1$, in which case the categories are pointed.
	
	\begin{lm}\label{lm:cyclic-pointed}
		Let $(\mathcal B,c)$ and $(\mathcal D,c')$ be pointed braided fusion categories whose groups of simple objects are cyclic of order $N$. Let $[X]$ and $[Y]$ be generators of these groups, respectively. Assume there exists $\lambda\in\C^\times$ such that
		\[
		c_{X,X}=\lambda\cdot \Id_{X\bt X},\qquad c'_{Y,Y}=\lambda\cdot \Id_{Y\bt Y}.
		\]
		Then there is a braided equivalence $F:\mathcal B\xrightarrow{\simeq}\mathcal D$ such that $F(X^{\bt i})\cong Y^{\bt i}$ for all $0\leq i<N$.
	\end{lm}
	\begin{proof}
		Identify both groups with $\Z/N$ by sending $[X]$ and $[Y]$ to $1$. Under these identifications, consider the quadratic forms associated to the braidings $c$ and $c'$, respectively: 
		\[
		\mathfrak q,\mathfrak q'\colon \Z/N\ra\C^\times. 
		\]
		By \cite[Section~8.4]{EGNO}, $\mathfrak q(1)$ is the scalar by which $c_{X,X}$ acts on $X\bt X$, and similarly for $\mathfrak q'(1)$ and $c'_{Y,Y}$. Hence, for each $0\leq i<N$, we have
		\[
		\mathfrak q(i)=\mathfrak q(1)^{i^2}=\lambda^{i^2}
		=\mathfrak q'(1)^{i^2}=\mathfrak q'(i).
		\]
		It follows from \cite[Theorem~8.4.12]{EGNO} that the group isomorphism sending $[X]$ to $[Y]$ induces a braided equivalence $F:\mathcal B\xrightarrow{\simeq}\mathcal D$ such that $F(X)\cong Y$. Thus $F(X^{\bt i})\cong F(X)^{\bt i}\cong Y^{\bt i}$ for $0\leq i<N$.
	\end{proof}
	
	For $k\geq2$, we first determine the Hecke parameter and the corresponding determinant scalar for the quantum group category.
	
	\begin{lm}\label{lm:QG-invariants}
		Let $N,k\geq2$ and $p=N+k$. Let $q,y\in\C^\times$ satisfy $\operatorname{ord}(q^2)=p$ and $y^N=q$, and let $Y=L_q(\om_1)$ in $(\CC(\sl_N,q),c^y)$. Then $Q=q^{-2}$ is a Hecke parameter of $(\CC(\sl_N,q),c^y)$ with respect to $Y$. For this choice of $Q$, the determinant scalar in \eqref{eq:KW-determinant} is
		\[
		\tau_{(\CC(\sl_N,q),c^y)}=s^{-N}=(-1)^Nq^{1-N},\quad \mathrm{where}\quad s=-y^{N-1}.
		\]
	\end{lm}
	\begin{proof}
		From the calculation on $Y\bt Y$ in the proof of Proposition~\ref{thm:KW-comparison}, we see that $Q=q^{-2}$ is a Hecke parameter of $(\CC(\sl_N,q),c^y)$ and $c^y_{Y,Y}=sT$, where $T$ is given by \eqref{eq:KW-generator}. Choose $\iota:\vac\ra Y^{\bt N}$ and $\pi:Y^{\bt N}\ra\vac$ such that $\pi\circ\iota=\Id_\vac$. By the hexagon identity, $c^y_{Y,Y^{\bt N}}=s^NT_N\cdots T_1$. Substituting into \eqref{eq:KW-determinant} and using naturality of $c^y$ with respect to $\iota$, we have
		\begin{align*}
			\tau_{(\CC(\sl_N,q),c^y)}\Id_Y
			&=s^{-N}(\pi\bt\Id_Y)\circ c^y_{Y,Y^{\bt N}}\circ(\Id_Y\bt\iota)\\
			&=s^{-N}(\pi\bt\Id_Y)\circ(\iota\bt\Id_Y)=s^{-N}\Id_Y.
		\end{align*}
		Since $y^N=q$, we have $s^{-N}=(-1)^Nq^{1-N}$.
	\end{proof}
	
	\begin{assu}\label{ass:recognition}
		Let $N\geq2$ and $k\geq1$ be integers, and let $p=k+N$. Let $\mathcal B$ be a ribbon fusion category with braiding $c$ and twist $\theta$. Assume that the following conditions hold.
		\begin{enumerate}[label=(\arabic*)]
			\item The parameters $q,y\in\C^\times$ satisfy \eqref{eq:order-q-square}
			\[
			\operatorname{ord}(q^2)=p,\qquad y^N=q.
			\]
			\item There exists an isomorphism of based rings
			\[
			\mathfrak F(\mathcal B)\xrightarrow{\simeq}\mathfrak F(L_k(\sl_N)),\qquad
			[X_\la]\mapsto[\hat L_k(\la)],\quad \la\in P^k_+,
			\]
			where $X_\la$ denotes the simple object of $\mathcal B$ corresponding to $\la$. Let $X=X_{\om_1}$. As before, $\om_2=0$ when $N=2$.
			\item Let $\lambda_-$ be the eigenvalue of $c_{X,X}$ on the summand $X_{\om_2}$ of $X\bt X$, in the sense of Definition~\ref{df:braiding-eigenvalue}. The twist of $X$ and this eigenvalue are given by
			\[
			\theta_X=y^{N^2-1},\qquad \lambda_-=-y^{-N-1}.
			\]
			If $k\geq2$, the eigenvalue of $c_{X,X}$ on $X_{2\om_1}$ is $\lambda_+=y^{N-1}$.
		\end{enumerate}
	\end{assu}
	
	\begin{thm}\label{thm:recognition}
		Under Assumption~\ref{ass:recognition}, there is a ribbon equivalence
		\[
		F:\mathcal B\xrightarrow{\simeq}(\CC(\sl_N,q),c^y)
		\]
		such that $F(X_\la)\cong L_q(\la)$ for every $\la\in P^k_+$.
	\end{thm}
	\begin{proof}
		Denote $\mathcal D=(\CC(\sl_N,q),c^y)$ and let $Y=L_q(\om_1)\in \mathcal D$. First we consider the case when $k=1$. By Assumption~\ref{ass:recognition} (2) and the fusion tensor identities \eqref{eq:fundamental-fusion-rule} and \eqref{eq:QG-fusion-rule}, we have
		\[
		X^{\bt i}\cong X_{\om_i},\qquad Y^{\bt i}\cong L_q(\om_i),\qquad 1\leq i<N,
		\]
		and the $N$-th tensor power is the unit in each category. Since $P^1_+=\{0,\om_1,\ldots,\om_{N-1}\}$, $\mathcal B$ and $\mathcal D$ are pointed, with cyclic groups of simple objects generated by $[X]$ and $[Y]$, respectively. Then by Assumption~\ref{ass:recognition} (3) and Proposition~\ref{prop:QG-eigenvalues} (2), we have 
		\[
		c_{X,X}=\lambda_-\cdot \Id_{X\bt X},\qquad c^y_{Y,Y}=\lambda_-\cdot \Id_{Y\bt Y},\qquad
		\lambda_-=-y^{-N-1}.
		\] 
		It follows from Lemma~\ref{lm:cyclic-pointed} that there is a braided equivalence $F:\mathcal B\xrightarrow{\simeq}\mathcal D$ such that $F(X_{\om_i})\cong L_q(\om_i)$ for $1\leq i<N$.
		
		Now assume $k\geq2$. Let $Q=q^{-2}$ and $s=-y^{N-1}$. By Lemma~\ref{lm:QG-invariants}, $Q$ is a Hecke parameter of $\mathcal D$ and the corresponding determinant scalar is $\tau_{\mathcal D}=s^{-N}$. To apply Proposition~\ref{thm:KW-comparison}, it remains to show that $Q$ is a Hecke parameter of $\mathcal B$ and that $\tau_{\mathcal B}=\tau_{\mathcal D}$.
		
		By Assumption~\ref{ass:recognition} (2) and \eqref{eq:fundamental-fusion-rule}, we have $X\bt X\cong X_{2\om_1}\oplus X_{\om_2}$. Let $a\in\End_{\mathcal B}(X\bt X)$ be the projection onto $X_{\om_2}$ and let $T=Qa-(\Id-a)$ as in \eqref{eq:KW-generator}. Then by Assumption~\ref{ass:recognition} (3), $\lambda_-=sQ$ and $\lambda_+=-s$. It follows that
		\[
		c_{X,X}=\lambda_-a+\lambda_+(\Id-a)=s\bigl(Qa-(\Id-a)\bigr)=sT.
		\]
In particular, each $T_i$ acts as $s^{-1}$ times the braiding on the $i$-th and $(i+1)$-st factors of $X^{\bt m}$. Hence $T_i$ satisfy the braid relations in \eqref{eq:KW-Hecke-relations}. Moreover, since $T=s^{-1} c_{X,X}$ acts as $Q \cdot \Id$ on $X_{\om_2}$ and as $-\Id $ on $X_{2\om_1}$, we have $(T-Q\Id)(T+\Id)=0$. Therefore, $Q$ is a Hecke parameter of $\mathcal B$ by Definition~\ref{df:KW-Hecke-parameter}.
		
Now let $\iota:\vac\ra X^{\bt N}$ and $\pi:X^{\bt N}\ra\vac$ satisfy $\pi\circ\iota=\Id_\vac$. Similar to the proof of Lemma~\ref{lm:QG-invariants}, it follows from the hexagon identity of $c$ that $T_N\cdots T_1=s^{-N}c_{X,X^{\bt N}}$. Hence
		\[
		\tau_{\mathcal B}\Id_X
		=s^{-N}(\pi\bt\Id_X)\circ(\iota\bt\Id_X)
		=s^{-N}\Id_X.
		\]
	This shows $\tau_{\mathcal B}=s^{-N}=\tau_{\mathcal D}$. i.e., $\mathcal B$ and $\mathcal D$ have the same determinant scalars given by $Q$. 
	
Therefore, by Proposition~\ref{thm:KW-comparison}, there is a tensor equivalence
		\[
		F:\mathcal B\xrightarrow{\simeq}\mathcal D,\qquad F(X_\la)\cong L_q(\la),\quad \la\in P^k_+.
		\]
		
It remains to show $F$ preserves braiding and twist on $\mathcal B$ and $\mathcal D$.
Let $c'$ be the braiding on $\mathcal D$ transported from $c$ along $F$. 
Recall that $X\bt X\cong X_{2\om_1}\op X_{\om_2}$ and $Y\bt Y\cong L_q(2\om_1)\op L_q(\om_2)$. 
Since $F$ preserves tensor with $F(X_{2\om_1})\cong L_q(2\om_1)$ and $F(X_{\om_2})\cong L_q(\om_2)$, then by Assumption~\ref{ass:recognition} (3), $c'_{Y,Y}$ acts by $\lambda_+$ on $L_q(2\om_1)$ and by $\lambda_-$ on $L_q(\om_2)$. Hence $c'_{Y,Y}=c^y_{Y,Y}$ by \eqref{eq:QG-eigenvalues}. Applying the hexagon identities successively, we have
		\[
		c'_{Y^{\bt m},Y^{\bt n}}=c^y_{Y^{\bt m},Y^{\bt n}},\qquad m,n\geq1.
		\]
By Lemma~\ref{lm:least-power}, every simple object in $\mathcal D$ occurs in a tensor power of $Y$. Hence the braidings $c'$ and $c^y$ agree on every pair of simple objects in $\mathcal D$. This shows $F$ preserves braiding.
	
	Finally, let $\theta$ be the given twist on $\mathcal D$ and let $\theta'$ be the twist on $\mathcal D$ transported along $F$. Since $F(X)\cong Y$, it follows from Assumption~\ref{ass:recognition} (3), \eqref{eq:QG-twist} and \eqref{eq:casimirs} that $\theta'_Y=y^{N^2-1}\Id_Y=\theta_Y$. Hence $\theta'=\theta$ by Lemma~\ref{lm:ribbon-structures}, since every simple object of $\mathcal D$ occurs in some $Y^{\bt m}$. Therefore, $F$ is a ribbon equivalence.
	\end{proof}
	
	
	\section{Ribbon equivalence between $\Oord$ and the quantum group category of type $A$ }\label{sec:main-section}
	Throughout this section, we let $\g=\sl_N$ with $N\geq2$, and let $\ell$ be an admissible level. Since $h^\vee=N$, \eqref{eq:kappa} becomes
	\begin{equation}\label{eq:typeAlabels}
		\kappa=\ell+N=\frac{p}{v},\qquad p,v\in\Z_{\geq1},\quad \gcd(p,v)=1,\quad p\geq N.
	\end{equation}
	Consider the following positive integer
	\begin{equation}\label{eq:levelk}
		k=p-N\in\Z_{\geq0}.
	\end{equation}
	Then $\mathrm{Pr}^{\ell}_{\Z}=P^k_+$ \eqref{eq:PrZ}, \cite[(2.5)]{AvE23}. Here $\ell$ is the admissible level for the VOA $L_\ell(\sl_N)$, while $k$ is the positive integral level of the comparison strongly rational affine VOA $L_k(\sl_N)$ \cite{DLM97}. 
	
	We prove the ribbon equivalence between $\Oord$ and $(\CC(\sl_N,q),c^y)$ by applying the recognition result of Section~\ref{sec:recognition}, using Lemma~\ref{lm:adm-ribbon-ADE} and the braiding formula of Section~\ref{sec:admissible}.
	
	\subsection{Braiding eigenvalues in type $A$}\label{sec:main}
	We first recall the fusion rules at the admissible level.
	\begin{lm}[{\cite[Corollary~7.4]{Cre19}}]\label{lm:adm}
		For $\la,\mu\in\mathrm{Pr}^{\ell}_{\Z}=P^k_+$, we have the following decomposition of $L_\ell(\sl_N)$-modules
		\[
		\Lh{\la}\bt\Lh{\mu}\cong\bigoplus_{\nu\in P^k_+}N_{\la,\mu}^\nu\,\Lh{\nu},
		\]
		where $N_{\la,\mu}^\nu$ are the fusion coefficients of $L_k(\sl_N)$ \cite{V88,FZ92,Bea96}. In particular, the assignment $[\Lh{\la}]\mapsto[\hat L_k(\la)]$ for $\la\in P^k_+$ gives an isomorphism of based rings
		\[
		\mathfrak F(\Oord)\xrightarrow{\simeq}\mathfrak F(L_k(\sl_N)).
		\] 
	\end{lm}
	
	Now consider the following roots of unities
	\begin{equation}\label{eq:adm-q-y}
		q=e^{\pi iv/p},\qquad y=e^{\pi iv/(Np)}.
	\end{equation}
	Then $y^N=q$ and $q^2$ has order $p$, since $\gcd(p,v)=1$.
	
	Recall that $\om_1,\ldots,\om_{N-1}$
	are the fundamental weights of $\sl_N$. When $N=2$, we let
	$\om_2=0$.	For $k\geq1$, let $X=\Lh{\om_1}\in\Oord$. Since $\ell+N=p/v$, it follows from Lemma~\ref{lm:adm-ribbon-ADE}, \eqref{eq:conformal-weight} and \eqref{eq:casimirs} that the twist $\theta_X$ in the category $\Oord$ is given as follows: 
	\begin{equation}\label{eq:adm-vector-twist}
		\theta_X=e^{2\pi i h_{\om_1}}
		=\exp\left(\frac{\pi i v}{p}(\om_1\mid\om_1+2\rho)\right)
		=e^{\pi i v(N^2-1)/(Np)}=y^{N^2-1}.
	\end{equation}
	We now compute the braiding eigenvalues appearing in Assumption~\ref{ass:recognition} (3). 
	
	\begin{prop}\label{prop:adm-invariants}
		Let $k\geq1$ and let $X=\Lh{\om_1}\in\Oord$. With $y$ as in
		\eqref{eq:adm-q-y}, we have 
		\begin{enumerate}[label=(\arabic*)]
			\item If $k\geq2$, then $X\bt X\cong\Lh{2\om_1}\oplus\Lh{\om_2}$. The eigenvalues of $c_{X,X}$ on $\Lh{2\om_1}$ and $\Lh{\om_2}$ are $\lambda_+=y^{N-1}$ and $\lambda_-=-y^{-N-1}$, respectively.
			\item If $k=1$, then $X\bt X\cong\Lh{\om_2}$ and the eigenvalue of $c_{X,X}$ on $\Lh{\om_2}$ is $\lambda_-=-y^{-N-1}$.
		\end{enumerate}
	\end{prop}
	\begin{proof}
		Let $E=L(\om_1)=\C^N$ as in \eqref{eq:Lomega1wedge}. Recall from Subsection~\ref{sec:QG} that
		\begin{equation}\label{eq:EotimesE}
			E\o_\C E=S^2E\oplus\bigwedge\nolimits^2 E,\qquad
			S^2E\cong L(2\om_1),\quad \bigwedge\nolimits^2 E\cong L(\om_2).
		\end{equation}
		By Lemma~\ref{lm:adm} and \eqref{eq:fundamental-fusion-rule}, the multiplicity of $\Lh{\nu}$ in $X\bt X$ is $N_{\om_1,\om_1}^{\nu}=c_{\om_1,\om_1}^{\nu}$ 	for all $\nu\in P^k_+$.  By definition  $c_{\om_1,\om_1}^{\nu}$ is the multiplicity of $L(\nu)$ in $E\o_\C E$ as an $\sl_N$-module. 	Hence $\Lh{\om_2}$ occurs with multiplicity one in $X\bt X$ for every $k\geq1$, in view of \eqref{eq:EotimesE}. 
		
		Since $(2\om_1\mid\theta)=2$, we have $2\om_1\in P^k_+$ if and only if $k\geq2$. It follows from \eqref{eq:EotimesE} that 
		\[
		X\bt X\cong
		\begin{cases}
			\Lh{2\om_1}\oplus\Lh{\om_2}, & k\geq2,\\
			\Lh{\om_2}, & k=1.
		\end{cases}
		\]
		
		We now determine the sign $\epsilon_Z$ in Corollary~\ref{coro:adm-braiding-eigenvalues}.
		By \eqref{eq:adm-T} and Lemma~\ref{thm:adm-bottom}, we have the $\g$-module isomorphisms
		\[
		T_{\om_1,\om_1}
		\xrightarrow[\sim]{\widetilde\varphi_{\om_1,\om_1}}
		X\odot X
		\xrightarrow[\sim]{\eta_{\om_1,\om_1}}
		\Om(X\bt X)\cong
		\begin{cases}
			L(2\om_1)\oplus L(\om_2), & k\geq2,\\
			L(\om_2), & k=1,
		\end{cases}
		\]
		where the last isomorphism follows from \eqref{eq:adm-Zhu-simples}.
		
		(1) Consider the case when $k\geq2$. By \eqref{eq:EotimesE}, we have $\dim T_{\om_1,\om_1}=\dim(E\o_\C E)=N^2$.
		Hence the quotient map $\pi_{\om_1,\om_1}:E\o_\C E\twoheadrightarrow T_{\om_1,\om_1}$ in \eqref{eq:Tlamu} is an isomorphism.
		Recall that  by \eqref{eq:Tlala}, the induced flip $\bar\tau$ acts by $1$ on the image of $S^2E\cong L(2\om_1)$ and by $-1$ on the image of $\bigwedge\nolimits^2 E\cong L(\om_2)$.
		Substituting these signs into $	\lambda_Z=\epsilon_Ze^{\pi i(h_\nu-2h_\la)}$ \eqref{eq:adm-braiding-eigenvalue}, with $Z=\Lh{2\om_1}$ and $Z=\Lh{\om_2}$, respectively, we have
		\begin{align*}
			\lambda_+&=e^{\pi i(h_{2\om_1}-2h_{\om_1})}
			=e^{\pi i v(N-1)/(Np)}=y^{N-1},\\
			\lambda_-&=-e^{\pi i(h_{\om_2}-2h_{\om_1})}
			=-e^{-\pi i v(N+1)/(Np)}=-y^{-N-1},
		\end{align*}
		where we used \eqref{eq:conformal-weight} and \eqref{eq:casimirs} for the explicit values of conformal weights. 
		
		(2) Consider the case when $k=1$. Then $T_{\om_1,\om_1}\cong \Om(X\bt X)\cong L(\om_2)$. In particular, by \eqref{eq:EotimesE},  $\pi_{\om_1,\om_1}:E\o_\C E\twoheadrightarrow T_{\om_1,\om_1}$ vanishes on $S^2E\cong L(2\om_1)$ and maps $\bigwedge\nolimits^2 E\cong L(\om_2)$ isomorphically onto $T_{\om_1,\om_1}$. Then by \eqref{eq:Tlala}, $\tau$ acts as $-1$ on $T_{\om_1,\om_1}$. Thus,
		\[
		\lambda_-=-e^{\pi i(h_{\om_2}-2h_{\om_1})}
		=-e^{-\pi i v(N+1)/(Np)}=-y^{-N-1},
		\]in view of \eqref{eq:adm-braiding-eigenvalue}. 
	\end{proof}
	
	\subsection{The main theorem and consequences}\label{sec:main}
	The proof of our main theorem uses the reconstruction theorem of Kazhdan and Wenzl \cite{KW93} for $k\geq2$ and the classification of pointed braided fusion categories for $k=1$. In both cases, we apply Theorem~\ref{thm:recognition} to $\mathcal B=\Oord$ with
	$X_\la=\Lh{\la}$ for each $\la\in P^k_+=\Pr^\ell_\Z$.
	
	\begin{thm}\label{thm:main}
		With 	$q=e^{\pi iv/p}$ and $ y=e^{\pi iv/(Np)}$ as in \eqref{eq:typeAlabels} and \eqref{eq:adm-q-y}, there is a ribbon equivalence of  categories
		\begin{equation}\label{eq:equivalenceF}
			F:\Oord\xrightarrow{\simeq}(\CC(\sl_N,q),c^y)
		\end{equation}
		such that $F(\Lh{\la})\cong L_q(\la)$ for every $\la\in P^k_+$.
	\end{thm}
	
	\begin{proof}[Proof of Theorem~\ref{thm:main}]
		If $k=0$, then $P^0_+=\{0\}$. Both categories have only the tensor unit as a simple object, and hence are ribbon equivalent to $\mathrm{Vec}$.
		Now assume $k\geq1$ and $\mathcal B=\Oord$. We show that all the conditions in Assumption~\ref{ass:recognition} are satisfied. 
		
		Indeed, since $\gcd(p,v)=1$, $q^2$ has order $p=k+N$, and $y^N=q$. Thus Assumption~\ref{ass:recognition} (1) holds. By Lemma~\ref{lm:adm-ribbon-ADE} applied to $\g=\sl_N$, $\Oord$ is a ribbon fusion category. Then the isomorphism of based rings of the fusion rings in Lemma~\ref{lm:adm} gives Assumption~\ref{ass:recognition} (2) with $\mathcal B=\Oord$ and $X_\la=\Lh{\la}$. 
		Finally, with the braiding eigenvalues $\la_\pm$ given by Proposition~\ref{prop:adm-invariants} and twist $\theta_X$ given by \eqref{eq:adm-vector-twist} , we see that the Assumption~\ref{ass:recognition} (3) holds. 
		The required ribbon equivalence $F$ now follows from Theorem~\ref{thm:recognition}.
	\end{proof}
	
	When $v=1$ \eqref{eq:kappa}, we have $p=k+N$ and $\ell=k$. Then, for  $k\geq1$, we have
	\[
	L_\ell(\sl_N)=L_k(\sl_N)
	\]
	is the strongly rational affine VOA at positive integral level. Its category of ordinary modules $\Oord$ is the semisimple category $\mathrm{KL}_k(\sl_N)$ consists of finite direct sums of integrable highest-weight $\widehat{\sl}_N$-modules of level $k$. 
	
	Theorem~\ref{thm:main} gives the following ribbon equivalence with the semisimplified category of tilting modules for $U_q(\sl_N)$ at a root of unity. This recovers the type~$A$ positive integral level case of the correspondence obtained by Kazhdan--Lusztig and Finkelberg \cite{KL93a,KL93b,KL94a,KL94b,F96,F13}.
	
	\begin{coro}
		Let $N\geq2$ and $k\geq1$ be integers, and let $q=e^{\pi i/(k+N)}$ and $y=e^{\pi i/(N(k+N))}$.
		Then there is an equivalence of ribbon categories
		\[
		F:\mathrm{KL}_k(\sl_N)\xrightarrow{\simeq}(\CC(\sl_N,q),c^y),\qquad F(\hat L_k(\la))\cong L_q(\la),\quad \la\in P^k_+. 
		\]
	\end{coro}
	
	We can also use Theorem~\ref{thm:main} to determine the modularity of $\Oord$ in type $A$.
	We use Rowell's formulation of Brugui\`eres' type~$A$ modularity criterion
	for quantum group module categories \cite[Section~5]{Bru00}.
	
	\begin{lm}[{\cite[Section~4.2]{Row06}}]\label{lm:qg-modularity}
		Let $p>N$, let $\gcd(v,p)=1$, and let $q=e^{\pi i v/p}$ and $y=e^{\pi i v/(Np)}$. Then the ribbon category $(\CC(\sl_N,q),c^y)$ is modular if and only if $\gcd(v,Np)=1$.
	\end{lm}
	
	\begin{coro}\label{coro:modularity}
		Under the assumptions of Theorem~\ref{thm:main}, the ribbon fusion category $\Oord$ is modular if and only if $\gcd(v,N)=1$.
	\end{coro}
	\begin{proof}
		If $\gcd(v,N)=1$, it was proved  by Creutzig that  $\Oord$ is modular \cite[Theorem~3]{Cre19}. 
		
		Conversely, suppose that $\Oord$ is modular. In view of  \eqref{eq:typeAlabels}, we must have $p\geq N$. 
		If $p=N$, then $\gcd(v,N)=1$ since $\gcd(v,p)=1$ \eqref{eq:typeAlabels}. If $p>N$, then the ribbon equivalence \eqref{eq:equivalenceF} implies that $(\CC(\sl_N,q),c^y)$ is modular. By Lemma~\ref{lm:qg-modularity}, we have $\gcd(v,Np)=1$, and hence $\gcd(v,N)=1$.
	\end{proof}
	
	In type $A$, \cite[Conjecture~1.2]{CHY18} states that $\Oord$ is modular whenever $N$ or $v$ is odd. However, this is false by Corollary~\ref{coro:modularity}.
	
	\begin{remark}
For braided tensor categories of type~$C$, Tuba and Wenzl proved a recognition theorem analogous to Theorem~\ref{thm:recognition} \cite[Theorem~9.4(a)]{TW05}. Hence to obtain a type~$C$ analogue of Theorem~\ref{thm:main}, one would need to verify the required fusion rules and ribbon structure for $\Oord$ and compare the twist and the self-braiding eigenvalues with the quantum group side. Our  braiding formula \eqref{eq:adm-braiding} applies to $\mathfrak{sp}_{2n}$ and gives the eigenvalues on the VOA side. 

However, the argument in \cite[Lemma~6.3(a)]{TW05} leaves out some admissible levels. Hence we leave type~$C$ and the other types for future work.
		
	\end{remark}

\end{document}